\documentclass[12pt]{article}
\usepackage[left=1in,top=1in,right=1in,bottom=1in]{geometry}
\usepackage{amssymb,amsfonts,amsmath,amsthm}
\usepackage{caption,graphicx,enumerate,xcolor}
\usepackage{hyperref,doi}
\numberwithin{equation}{section}
\numberwithin{figure}{section}
\numberwithin{table}{section}

\newcounter{labelflag} 

\newcommand{\Label}[1]{
                       \ifnum\thelabelflag=1 
                          \ifmmode  
                             \makebox[0in][l]{\qquad\fbox{\rm#1}}
                          \else
                             \marginpar{\vspace{0.7\baselineskip}
                                        \hspace{-1.1\textwidth}
                                        \fbox{\rm#1}}
                          \fi 
                       \fi
                       \label{#1} 
                      }

\newcommand{\eps}{\varepsilon}
\newcommand{\CO}{\mbox{$\mathrm{CO}_2$}}

\newcounter{remarkcounter}
\theoremstyle{definition}
\newtheorem{remark}[remarkcounter]{Remark}
\numberwithin{remarkcounter}{section}

\title{\vspace{-2em}\Large{Dynamical systems analysis of the Maasch--Saltzman\\ model for glacial cycles}}

\author{\normalsize
Hans Engler\footnotemark[1], 
Hans G. Kaper\footnotemark[1],
Tasso J. Kaper,\footnotemark[2] and
Theodore Vo\footnotemark[2]
}

\date{\vspace{-0.5em}\small{\today}}

\begin{document}
\maketitle
\bibliographystyle{abbrv}

\renewcommand{\thefootnote}{\fnsymbol{footnote}}
\footnotetext[1]{Department of Mathematics and Statistics, Georgetown University, Washington, DC 20057, and Mathematics and Climate Research Network (MCRN, \url{https://mathclimate.org})} 
\footnotetext[2]{Department of Mathematics and Statistics, Boston University, Boston, MA 02215}
\renewcommand{\thefootnote}{\arabic{footnote}}

\vspace{-1em}
\begin{abstract}
\noindent
This article is concerned with the internal dynamics of a conceptual model proposed by Maasch and Saltzman [J. Geophys. Res., $\bf 95, D2$ (1990) 1955-1963] to explain central features of the glacial cycles observed in the climate record of the Pleistocene Epoch. It is shown that, in most parameter regimes, the long-term system dynamics occur on certain intrinsic two-dimensional invariant manifolds in the three-dimensional state space. 
These invariant manifolds are slow manifolds when the characteristic time scales for the total global ice mass and the volume of North Atlantic Deep Water are well-separated, and they are center manifolds when the characteristic time scales for the total global ice mass and the volume of North Atlantic Deep Water are comparable. 
In both cases, the reduced dynamics on these manifolds are governed by Bogdanov-Takens singularities, and the bifurcation curves associated to these singularities organize the parameter regions in which the model exhibits glacial cycles.
\end{abstract}

\section{Introduction} 		\label{s-Introduction}
The dynamics of glacial cycles during the Pleistocene Epoch---the period from approximately 2.6~million years before present (2.6\,Myr~BP) until approximately 11.7~thousand years before present (11.7\,Kyr~BP)---are of great current interest in the geosciences community, see \cite{Saltzman2001}, \cite[\S 11]{Dijkstra2013}  and \cite[\S 12.3]{MarshallPlumb2007}. The geological record shows cycles of advancing and retreating continental glaciers, mostly at high latitudes and high altitudes, and especially in the Northern Hemisphere. The typical temperature pattern inferred from proxy data resembles that of a sawtooth wave, where a slow glaciation is followed by a rapid deglaciation. In the early Pleistocene (until approximately 1.2\,Myr~BP), the period of a glacial cycle averaged 40\,Kyr; after a transition period of approximately 400\,Kyr, the glacial cycles had a noticeably greater amplitude and their period averaged 100\,Kyr. Although the periods appear to correlate to the cycles of the orbital forcing (Milankovitch theory \cite{Milankovitch1941}), the evidence is subject to debate~\cite[\S~11.8]{KaperEngler2014}, and there is currently no widely-accepted explanation for the mid-Pleistocene transition, when the period of the cycles changed from 40\,Kyr to 100\,Kyr. Several models have been proposed to explain the various observations; see, for example, \cite{Ashkenazy2004, Ashwin2015, Clark1999, Ghil1994, Huybers2005,Huybers2006,Huybers2007, MaaschSaltzman1990, Paillard1995, Paillard1998, Paillard2001, PaillardParrenin2004, Raymo1997, Saltzman1987, SaltzmanMaasch1988, SaltzmanMaasch1991, Shackleton2000, Tziperman2003}. We refer the reader to~\cite{Crucifix2012} for an overview of these various modeling efforts and to \cite{Saltzman2001} for a general introduction to paleoclimate modeling. 
The present investigation focuses on the internal dynamics of the conceptual model developed by Maasch and Saltzman~\cite{MaaschSaltzman1990,SaltzmanMaasch1988}. 

\subsection{The Maasch--Saltzman Model\label{ss-MSModel}}
The Maasch and Saltzman (MS) model is based on physical arguments and emphasizes the role 
of atmospheric \CO\ 
in the development and evolution 
of the glacial cycles.
In nondimensional form,
it consists of the following 
three ordinary differential equations:
\begin{equation}
\begin{split}
\dot{x} &= - x - y , \\
\dot{y} &= r y - p z + s z^2 - y z^2 , \\
\dot{z} &= - qx - q z .
\end{split}
\Label{MS-xyz}
\end{equation}
The state variables 
$x$, $y$, and $z$ represent
the anomalies 
(deviations from long-term averages) 
of the total global ice mass,
the atmospheric \CO\ concentration, and
the volume 
of the North Atlantic Deep Water (NADW),
respectively.
The latter 
is a measure of the strength 
of the North Atlantic overturning circulation 
and thus
of the strength of the oceanic \CO\ pump.
The parameters $p$, $q$, $r$, and $s$ are
combinations of various physical parameters.
They are all positive, 
with $q > 1$ 
for physical reasons.
The derivation 
of the model
is given 
in~\cite[\S 2]{SaltzmanMaasch1988}.

In \cite{MaaschSaltzman1990}, Maasch and Saltzman showed computationally that the model \eqref{MS-xyz} exhibits oscillatory behavior with dominating periods of 40~Kyr in response to insolation forcing with such periods, and limit cycles with 100~Kyr periods if $p=1, q=1.2, r=0.8,$ and $s=0.8$ in the absence of forcing.
They also showed that a transition from a 40 Kyr cycle to a 100 Kyr cycle can be achieved by slowly varying the parameters $p$ and $r$ across a certain threshold.

\begin{figure}[h!]
\centering
\includegraphics[width=5in]{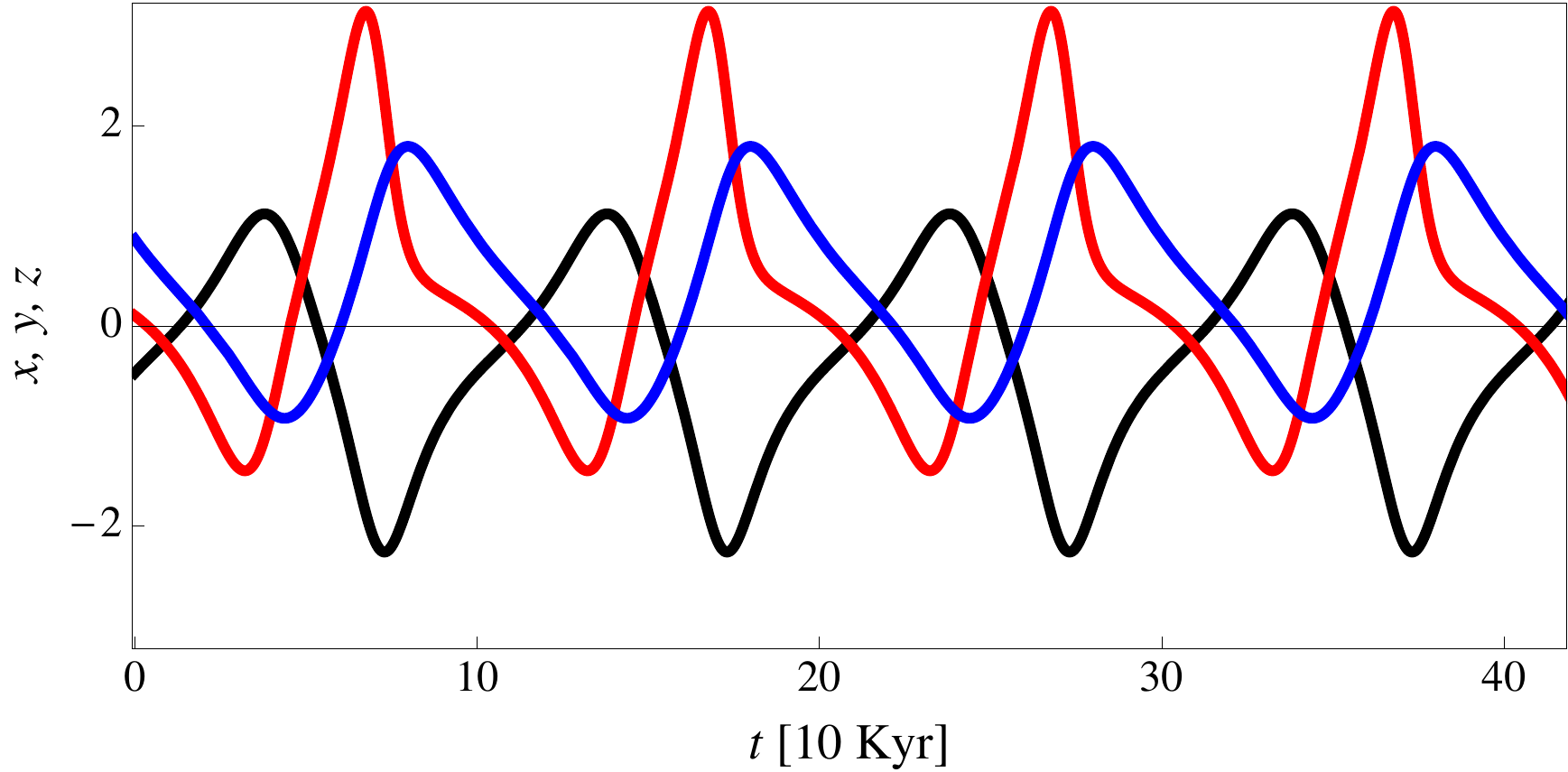}
\caption{Limit cycle of \eqref{MS-xyz} at $p =1.0$, $q = 1.2$, $r = 0.8$, $s = 0.8$.
The three curves represent the total ice mass (black), atmospheric \CO\ concentration (red), and volume of NADW (blue).
\label{f-Intro-LimitCycle}}
\end{figure}

Figure \ref{f-Intro-LimitCycle} shows a representative 100 Kyr limit cycle. 
Each cycle is clearly asymmetric: a rapid deglaciation is followed by a slow glaciation. 
This asymmetry arises in \eqref{MS-xyz} for $s>0$.
Also, the three variables are properly correlated:
as the concentration of the atmospheric \CO\ (a greenhouse gas)  increases, the climate gets warmer, and the total ice mass decreases (deglaciation); as the volume of NADW increases, the strength of the North Atlantic overturning circulation increases, more atmospheric \CO\ is absorbed by the ocean and, consequently, the atmospheric \CO\ concentration decreases. 

In this article, we present a dynamical systems analysis of the internal dynamics of the Maasch--Saltzman (MS) model \eqref{MS-xyz}. We identify the Bogdanov-Takens (BT) points \cite{Bogdanov1975,Broer2001,Guckenheimer1983,Kuznetsov2004,Takens1974} that act as organizing centers in the parameter space for all of the equilibria, limit cycles, homoclinic orbits, and their bifurcations. In addition to being of intrinsic interest, our analysis of the internal dynamics of \eqref{MS-xyz} will be instrumental for investigations of the effects of time-dependent forcing, especially of orbital (Milankovitch) forcing, and of the effects of slowly-varying parameters $p$ and $r$.

Two observations are useful for the analysis of the MS model. First, in the special case $s=0$, the system \eqref{MS-xyz} reduces to 
\begin{equation}
\begin{split}
\dot{x} &= - x - y , \\
\dot{y} &= r y - p z - y z^2 , \\
\dot{z} &= - q x - q z,
\end{split}
\Label{MSsym-xyz-1}
\end{equation}
which possesses a reflection symmetry; if $(x, y, z)$ is a solution then so is $(-x, -y, -z)$.
Hence, it will be useful to study the {\it symmetric MS} model~\eqref{MSsym-xyz-1} and to use the results
to understand how the dynamics change 
for $s>0$, as the physical symmetry is broken.

Second, 
the parameter $q$ 
is essentially the ratio
of the characteristic time scales 
for the total global ice mass~($x$)
and the volume of NADW~($z$).
It turns out to be useful 
to consider first
the asymptotic case 
where $q$ is finite
but large, 
denoted by $q \gg 1$.
In this case, 
the MS model
(either in its original form 
or in its symmetric form)
becomes a {\it slow--fast system}, 
where the fast variable~$z$
is slaved 
to the slow variables~$x$ and~$y$.
The insights gained
from this analysis
will then be a useful guide
for understanding the dynamics
for all finite $q>1$.

\subsection{Summary of the Results\label{ss-SummaryofResults}}
The first results are for the slow--fast regime of the symmetric MS model \eqref{MSsym-xyz-1}. In this regime, the system is (2+1)-dimensional, with two slow variables and one fast variable. 
We refer to this as the slow--fast symmetric model.
We show that there is a family of two-dimensional slow invariant manifolds to which all solutions quickly relax, and we study the dynamics on the slow manifolds. The central feature is a $\mathbb{Z}_2$-symmetric BT bifurcation point, from which all bifurcation curves emanate. The curves of Hopf bifurcations, homoclinic bifurcations, and saddle-node bifurcations of limit cycles determine the regions in parameter space where the stable limit cycles exist. In addition, since all solutions relax quickly to the slow manifolds, one can determine the basins of attraction of the various limit cycles. These first results build naturally on the recent analysis of the symmetric MS model \eqref{MSsym-xyz-1} in the limit $q = \infty$ \cite{EKKV2017}. 

The second results concern the effects of asymmetry ($s>0$). In the regime of finite but large values of $q$ ($q \gg 1$), the system~\eqref{MS-xyz} with $s>0$ is also a slow--fast system. 
We refer to it as the slow--fast asymmetric system.
There is again a family of exponentially attracting two-dimensional, invariant slow manifolds, but the symmetry-breaking makes the dynamics on the slow manifolds more complex. With $s>0$, the limit cycles observed in \eqref{MS-xyz} are asymmetric, exhibiting a relatively rapid deglaciation and a relatively slow glaciation, as shown in Figure \ref{f-Intro-LimitCycle}.

With these results in hand, we are then in a position to analyze and visualize the dynamics of the full, asymmetric ($s>0$) MS model~\eqref{MS-xyz} for all $q>1$. 
We show that for all $q>1$ the system possesses a family of two-dimensional center manifolds toward which solutions relax. Moreover, the solutions of the full system may be accurately approximated by those of the reduced systems on the center manifolds for all $q>1$, and the manifold is at least $C^1$-smooth for all $q$ greater than a critical value $q_c(p,r,s)$. 
On the center manifolds, the system has a pair of BT singularities, and the bifurcation curves emanating from them organize the system dynamics, including the boundaries of the regions where stable limit cycles exist. 

\subsection{Outline of the Article\label{ss-Outline}}
The article is organized as follows. 
In Section~\ref{s-FSsym}, we present the analysis of the slow--fast symmetric system.
In Section~\ref{s-FSasym}, we present the analysis of the slow--fast asymmetric system.
Section~\ref{s-MS} presents the analysis of the full three-dimensional MS model. 
We conclude with a discussion in Section~\ref{s-Discussion}.
Appendix~\ref{sss-Q2homoclinics} provides details of the unfolding and Melnikov analysis for the persistence of homoclinic orbits in the slow--fast asymmetric system. 
Appendix~\ref{app:CM} presents essential information about the center manifolds in \eqref{MS-xyz}.

\section{Slow--Fast Dynamics of the Symmetric Model \label{s-FSsym}}
In this section, we analyze the symmetric MS model~\eqref{MSsym-xyz-1} for large $q$ ($q \gg 1$). The system is readily formulated as a slow--fast system, 
\begin{equation}
\begin{split}
\dot{x} &= - x - y , \\
\dot{y} &= r y - p z - y z^2 , \\
\eps \dot{z} &= - x - z ,
\end{split}
\Label{FSsym-xyz}
\end{equation}
where $\eps = 1/q$ is the small parameter, which measures the separation of time scales. We refer to~\eqref{FSsym-xyz} as the {\it slow--fast symmetric model}. Here, $x$ and $y$ are slow, and $z$ is fast. We assume that $p$ and $r$ are $\mathcal{O}(1)$ with respect to $\eps$.

\subsection{Slow Manifolds\label{ss-FSsymFSsym-Meps}}
With $\eps=0$, the system~\eqref{FSsym-xyz} has a critical manifold $\mathcal{M}_0 = \{ (x, y, z): z = - x\}$, which is invariant under the flow. Since $(\partial / \partial z) (-x-z) = -1$ for all $(x, y)$, $\mathcal{M}_0$ is normally attracting.
By Fenichel's Geometric Singular Perturbation Theory \cite{F1979,J1994,K1999}, normal hyperbolicity of $\mathcal{M}_0$ implies that, for any sufficiently small and positive $\eps$, there exists a family of persistent normally attracting invariant slow manifolds,
\begin{equation} 
\mathcal{M}_\eps = \{ (x, y, z): z = h_\eps (x, y) \} .
\end{equation}
The functions $h_\eps$ are $C^k$ for any $k > 0$. They are solutions of the invariance equation
\begin{equation}
\eps \frac{d}{dt} {h}_\eps (x, y) = - x - h_\eps (x, y),
\Label{FSsym-invariance}
\end{equation}
which satisfy $\displaystyle \lim_{\eps \to 0} h_{\eps}(x,y) = -x$.
The functions $h_\eps$ 
are identical to all orders 
in powers of~$\eps$
and differ only by terms of
$\mathcal{O} (e^{-c/\eps})$ 
as $\eps \to 0$, 
for some $c > 0$.
The expansion may be represented by 
\begin{equation}
h_\eps (x, y) = h_0 (x, y) + \eps h_1 (x, y) + \eps^2 h_2 (x, y) + \cdots .
\Label{FSsym-heps}
\end{equation}
The terms
in this expansion 
are found 
by substituting~(\ref{FSsym-heps})
into~(\ref{FSsym-invariance}) 
and equating coefficients 
of like powers of~$\eps$.
All terms in $h_\eps$
are of odd degree in the variables,
due to the $\mathbb{Z}_2$ symmetry
of~(\ref{FSsym-xyz}). 
The first few terms are
\begin{equation}
\begin{split}
h_0 (x, y) &= - x , \\
h_1 (x, y) &= - (x + y) , \\
h_2 (x, y) &=  - (x + y) + (ry + px- x^2 y) , \\
h_3 (x, y) &= - (1 - 2p + 4xy) (x + y) + (1 - r + x^2) (ry + px - x^2 y ) .
\end{split}
\Label{FSsym-hi}
\end{equation}

On the slow manifolds $\mathcal{M}_\eps$, the system~\eqref{FSsym-xyz} reduces to the planar system,
\begin{equation}
\begin{split}
\dot{x} &= - x - y , \\
\dot{y} &= r y - p h_\eps (x, y) - (h_\eps (x, y))^2 y ,
\end{split}
\Label{FSsym-xyheps}
\end{equation}
where $h_\eps$ is given by~\eqref{FSsym-heps} and \eqref{FSsym-hi}. The system \eqref{FSsym-xyheps} is the object of study in this section.

For completeness, we observe that the fast dynamics, along which solutions relax to $\mathcal{M}_\eps$, may be analyzed by introducing the fast time $\tau = t/\eps$ and rewriting \eqref{FSsym-xyz} as a fast system
\begin{equation}
\begin{split}
x' &= - \eps (x + y) , \\
y' &= \eps (r y - p z - y z^2) , \\
z' &= - x - z,
\end{split}
\Label{FSsym-xyz'}
\end{equation}
where the prime denotes differentiation with respect to $\tau$. Systems~\eqref{FSsym-xyz} and \eqref{FSsym-xyz'} are equivalent for all $\eps \not= 0$. In the limit as $\eps \to 0$, the fast system reduces to a single equation for $z$, with $x$ and $y$ constant in time. Let $(x_0,y_0,z_0)$ denote an arbitrary initial condition. In the fast time, $\tau$, the solution with initial condition $(x_0,y_0,z_0)$ relaxes to the point $(x_0,y_0,-x_0) \in \mathcal{M}_0$. Then, for $\eps$ small and positive, the $z$ component again relaxes quickly on the $\tau$ time scale, now to $\mathcal{M}_\eps$, and the $x$ and $y$ components will only change slowly guided by the dynamics on $\mathcal{M}_\eps$, see Section~\ref{ss-FSsym-Attraction}.

\subsection{Global Bifurcations\label{ss-FSsymBT}}
The trivial state, $P_0 = (0, 0)$,
is an equilibrium 
of the system~(\ref{FSsym-xyheps}) 
on $\mathcal{M}_\eps$
for all values of the parameters $(p,r)$.
In the regime $r > p$, 
there are also
equilibria at
$P_1 = (x_1^*, - x_1^*)$ and $P_2 = (x_2^*, - x_2^*)$
on $\mathcal{M}_\eps$,
where $x_1^* = \sqrt{r - p}$ 
and $x_2^* = - \sqrt{r - p}$.
The equilibria $P_1$ and $P_2$
emerge in a symmetric pitchfork
bifurcation from $P_0$
along the diagonal $r=p$.

\begin{figure}[h!]
\centering
\includegraphics[width=4in]{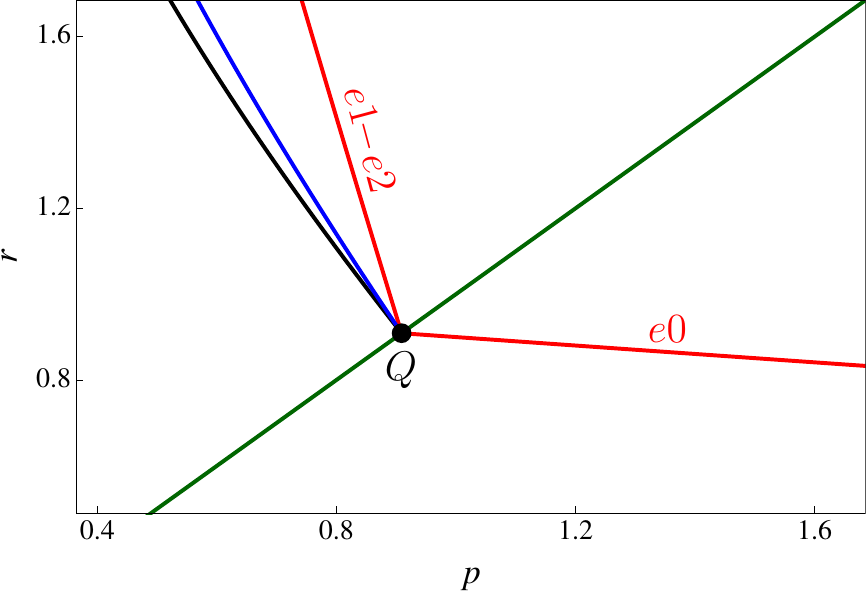}
\caption{Bifurcation structure of \eqref{FSsym-xyz} for $\eps=0.1$ ($q = 10$).}
\label{f-FSsym-Bifurcation}
\end{figure}

At the point $(p,r)=(\frac{1}{1+\eps}, \frac{1}{1+\eps})$ in the parameter space, the equilibrium $P_0$ of (\ref{FSsym-xyheps}) undergoes a $\mathbb{Z}_2$-symmetric BT bifurcation, since the Jacobian has a zero eigenvalue of geometric multiplicity two. This point is referred to as an organizing center \cite{Holmes1980}, and we denote it by $Q$. All bifurcation curves emanate from this point, including a pair of Hopf bifurcation curves, a homoclinic bifurcation curve, and a curve of saddle-node bifurcations of limit cycles; see \cite[\S8.4]{Kuznetsov2004} and \cite{Kuznetsov2005}. The organizing center 
and these bifurcation curves are shown in Figure~\ref{f-FSsym-Bifurcation}.
The representations of the global bifurcation curves emanating from the organizing centers, shown here, as well as those shown throughout this article, are obtained from numerical continuation using the software package AUTO~\cite{D1981,DCFKOPSWZ2007,DKK1991}.

The equilibrium $P_0$ undergoes a supercritical Hopf bifurcation along the (red) curve, which emanates to the right from the organizing center. The curve is given by
\begin{equation} 
{\rm e0} = \{ r = 1 - \eps p + \eps^3 p(p-1)+\mathcal{O}(\eps^4), \,  p > 1 - \eps + \eps^2 - \eps^3 + \mathcal{O}(\eps^4) \}.
\end{equation}
This curve is quadratic in $p$ and concave up, and the tangent line has slope $-\eps$ at the organizing center to leading order.
Similarly, $P_1$ and $P_2$ undergo subcritical Hopf bifurcations along the (red) curve emanating to the left from $Q$. The curve is given by
\begin{equation}
\mbox{ e1-e2} = \{ p = 1 + \eps (1 - 2r) + \eps^2 (1 - 2r) + \eps^3 (3 - 8r + 4r^2) + \mathcal{O}(\eps^4),\, r > 1 - \eps + \eps^2 - \eps^3 + \mathcal{O}(\eps^4)\}.
\end{equation}
This curve is quadratic in $r$ 
and concave up,
and the tangent line 
has slope $-1 / (2\eps)$
at the organizing center to leading order.

The homoclinic bifurcation curve
is shown in blue.
Its existence
is established 
using standard BT theory,
which entails an unfolding procedure 
followed by a Melnikov analysis.
This analysis shows that, to leading order, at the organizing center, the tangent line to the homoclinic bifurcation curve is
$r - (1-\eps) = -4 ( p - (1-\eps) )$, see also Section \ref{ss-BTanalysis}.

The curve of saddle-node bifurcations of limit cycles is shown in black. The existence of this second global bifurcation curve is a consequence of the $\mathbb{Z}_2$ symmetry and also follows from BT theory \cite{Kuznetsov2004,Kuznetsov2005}. The analysis also shows that, to leading order, this curve is tangent to the line $r - (1-\eps) \approx -3.03 ( p - (1-\eps))$ at the organizing center.

\begin{remark}
In the limit $\eps \to 0$ ($q \to \infty$), the slow manifolds $\mathcal{M}_\eps$ approach the critical manifold $\mathcal{M}_0$, and \eqref{FSsym-xyheps} reduces to the two-dimensional symmetric model studied in \cite[\S~4]{EKKV2017}.
\end{remark}

\begin{figure}[h!]
\centering
\includegraphics[width=4in]{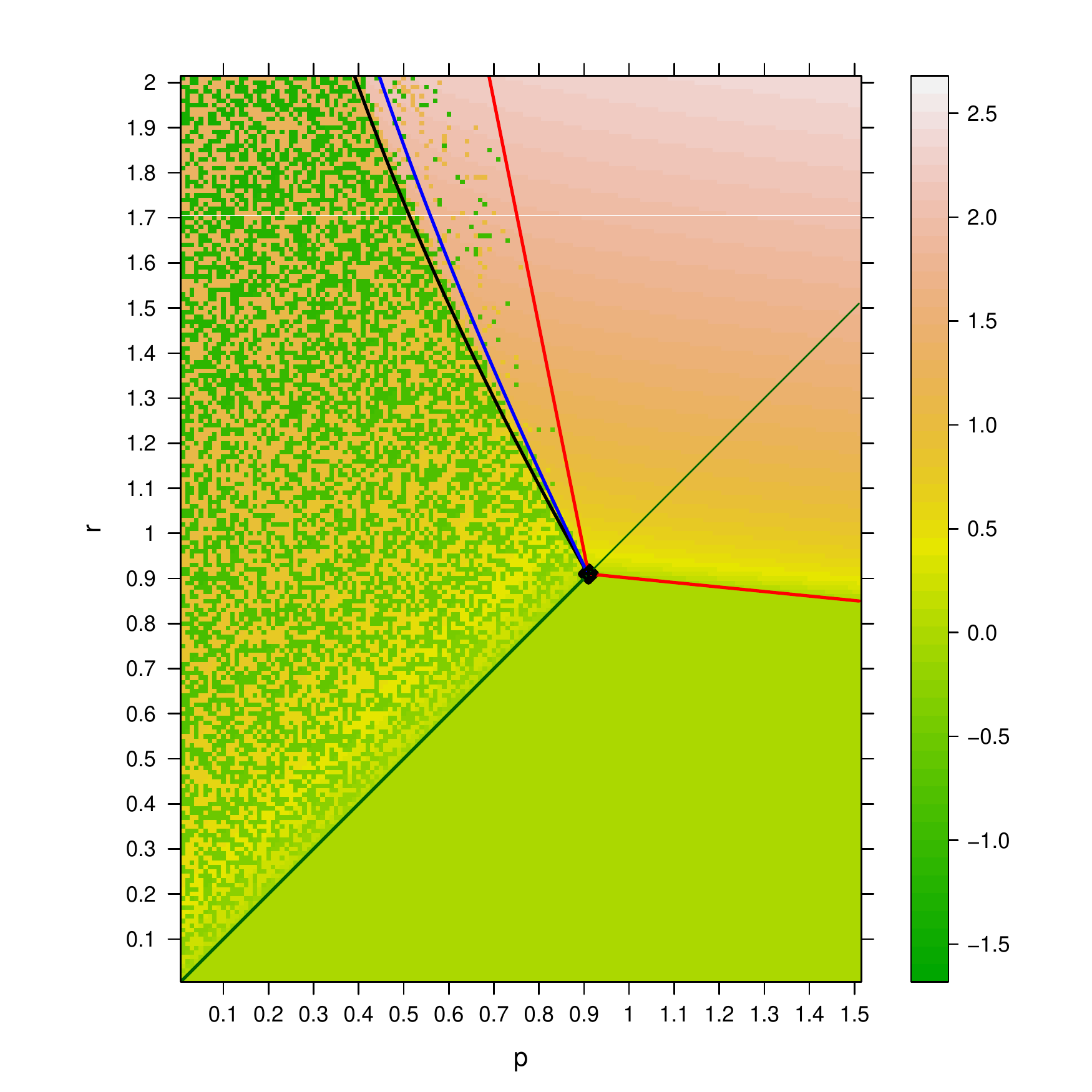}
\caption{Color map of $\overline{x} (p,r)$ for the system~(\ref{FSsym-xyz}) for $q = 10.0$ ($\eps=0.1$), where $\overline{x}$ is the function defined in~(\ref{xbar}), along with the bifurcation curves of the reduced fast system~(\ref{FSsym-xyheps}) with $h_\eps$ calculated up to and including $\mathcal{O}(\eps^3)$. \label{f-FSsym-HeatPlot}}
\end{figure}

Figure~\ref{f-FSsym-HeatPlot} shows the dynamics of \eqref{FSsym-xyz} in terms of the function $(p, r) \mapsto \overline{x} (p,r)$ for $(p, r) \in (0, 1.5) \times (0, 2)$ as a color map. Here, $\overline{x}$ is the function defined by 
\begin{equation}
\overline{x} = \lim \sup_{t \to \infty} x(t) ,
\Label{xbar}
\end{equation}
and the initial conditions were chosen randomly for each $(p,r)$.
Below the diagonal and e0, the color map is light green since all solutions approach $P_0$. As one crosses e0 with increasing $r$, the color changes from light green to orange and then to pink, indicating the presence of periodic orbits which are created in the supercritical Hopf bifurcation. Between the Hopf bifurcation curve e1-e2 and the homoclinic bifurcation curve, there is a similar shift to pink as $r$ increases. One also sees some green patches, indicating that some of the randomly chosen initial conditions lie in the basin of attraction of the stable equilibrium $P_2$. 
Next, in the region between the homoclinic bifurcation curve and the curve of saddle-node bifurcations of limit cycles, the color map has largely the same color, because solutions with initial conditions that lie inside the large unstable limit cycle approach one of the stable equilibria (green or orange), and those with initial conditions outside the unstable limit cycle approach the large stable limit cycle (pink). Finally, below the curve of saddle-node bifurcations of limit cycles, the color map consists entirely of green and orange, indicating that all of the solutions are attracted either to $P_1$ or to $P_2$, as expected since there are no stable limit cycles.

\subsection{Basins of Attraction\label{ss-FSsym-Attraction}}
An important feature of slow--fast systems like~(\ref{FSsym-xyz}) is that Geometric Singular Perturbation Theory provides insight into the dynamics not only of solutions on the slow manifold~$\mathcal{M}_\eps$, but also of solutions in a neighborhood of $\mathcal{M}_\eps$ and thus provides a way to explore the basins of attraction
of different solutions---such as equilibria and limit cycles---on $\mathcal{M}_\eps$.

In a neighborhood of $\mathcal{M}_\eps$, any solution $X(t) = (x(t), y(t), z(t))$ of~(\ref{FSsym-xyz}) decomposes into fast and slow components. The fast component is directed along the fast fiber $\mathcal{F}_\eps (b(t))$, which stands above the base point $b (t)$ of the fiber on $\mathcal{M}_\eps$. It decays exponentially fast toward the manifold $\mathcal{M}_\eps$, $\| X(t) - b(t) \| \le e^{-c t / \eps} \| X(0) - b(0) \|$ for some positive constant $c$, which is of $\mathcal{O} (1)$ as $\eps \to 0$. 
The family of fast fibers is invariant in the sense that $\mathcal{F}_\eps (b(t)) = \phi_t \left( \mathcal{F}_\eps (b(0)) \right)$ for all $t$, where $\phi_t$ is the time-$t$ flow map of~(\ref{FSsym-xyz}). The slow component of $X(t)$ lies in the tangent plane to the slow manifold~$\mathcal{M}_\eps$ at the base point~$b(t)$.

\begin{figure}[h!]
\centering
\includegraphics[width=3in]{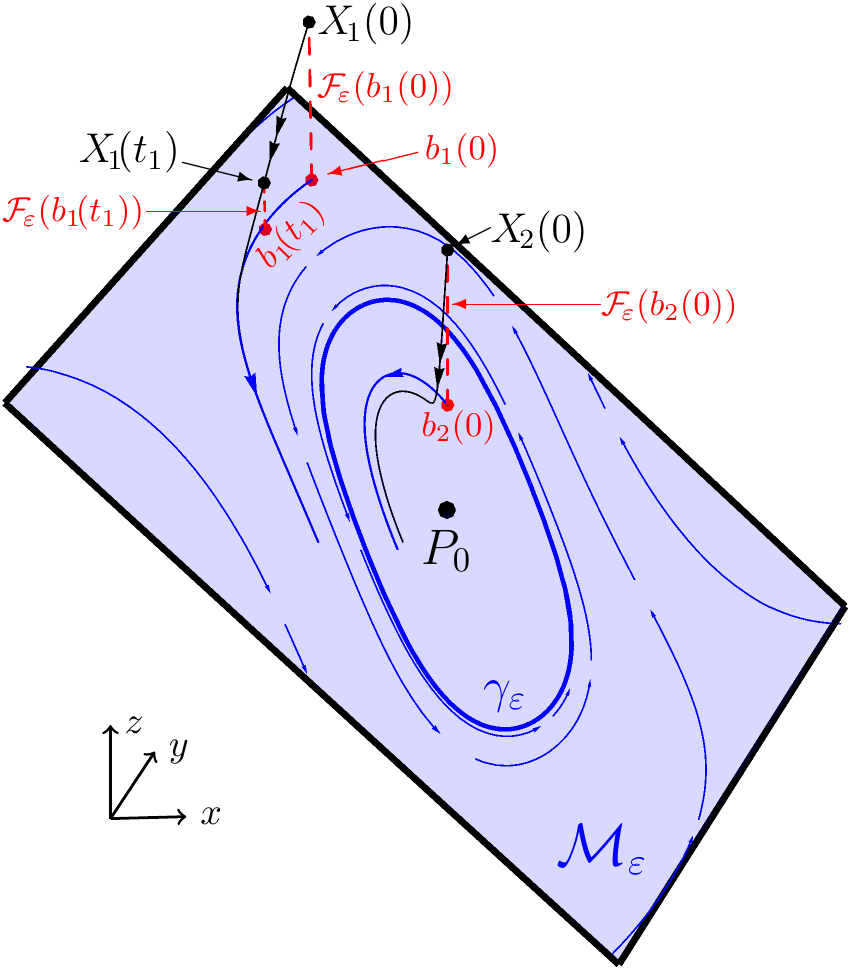}
\caption{
Slow--fast decomposition of typical solutions of (\ref{FSsym-xyheps})  in the basin of attraction $\mathcal{B}$ of the stable limit cycle $\gamma_\eps$ on $\mathcal{M}_\eps$, for parameters $(p,r)$ in the region between e0 and the diagonal where $P_0$ is an unstable spiral. The solution $X_1(t)$, with initial condition $X_1(0)$ on the fast stable fiber $\mathcal{F}_\eps(b_1(0))$, decomposes into a fast component that decays along the invariant family of fibers $\mathcal{F}_\eps(b_1(t))$ and a slow component that moves with the base point~$b_1(t)$. \label{f-FSsym-BasinRegionI}}
\end{figure}

We now demonstrate ---in two cases--- how one can use this decomposition to investigate the basins of attraction of solutions on $\mathcal{M}_\eps$. The first case is shown in Figure~\ref{f-FSsym-BasinRegionI}. Here, $(p,r)$ is in the region between e0 and the diagonal where $P_0$ is an unstable spiral and there is a unique stable
limit cycle $\gamma_\eps$ of~(\ref{FSsym-xyz}) around~$P_0$.

Let $\mathcal{B}$ denote the basin of attraction of $\gamma_\eps$.
The set $\mathcal{B} \vert_{\mathcal{M}_\eps}$ of all initial conditions
on $\mathcal{M}_\eps$ that lie in $\mathcal{B}$ is completely
determined by the analysis of the slow flow~(\ref{FSsym-xyheps})
on~$\mathcal{M}_\eps$.
For any initial condition $X(0) = (x(0), y(0), z(0))$
that lies near $M_\eps$ but not on it, there is a unique
fast stable fiber that contains $X(0)$.
Let $b(0)$ be the base point of this fiber on $\mathcal{M}_\eps$,
and denote the fiber by $\mathcal{F}_\eps (b(0))$.
Since any solution $X(t)$ near $\mathcal{M}_\eps$
decomposes into a fast component and a slow component,
we know immediately that $X(0)$ lies in $\mathcal{B}$
whenever $b(0)$ lies in $\mathcal{B} \vert_{\mathcal{M}_\eps}$.
Moreover, since the fast stable fibers completely foliate the
neighborhood of $\mathcal{M}_\eps$ 
and each initial condition $X(0)$ near $\mathcal{M}_\eps$
lies on a unique fiber, one may apply the above analysis
to each initial condition $X(0)$ near $\mathcal{M}_\eps$.
Therefore, the basin of attraction $\mathcal{B}$ of $\gamma_\eps$
contains all solutions near $\mathcal{M}_\eps$,
except those that start exactly on the stable fiber $\mathcal{F}_\eps (P_0)$
with base point $P_0$.
($P_0$ is not in $\mathcal{B} \vert_{\mathcal{M}_\eps}$.)

\begin{figure}[h!]
\centering
\includegraphics[width=3in]{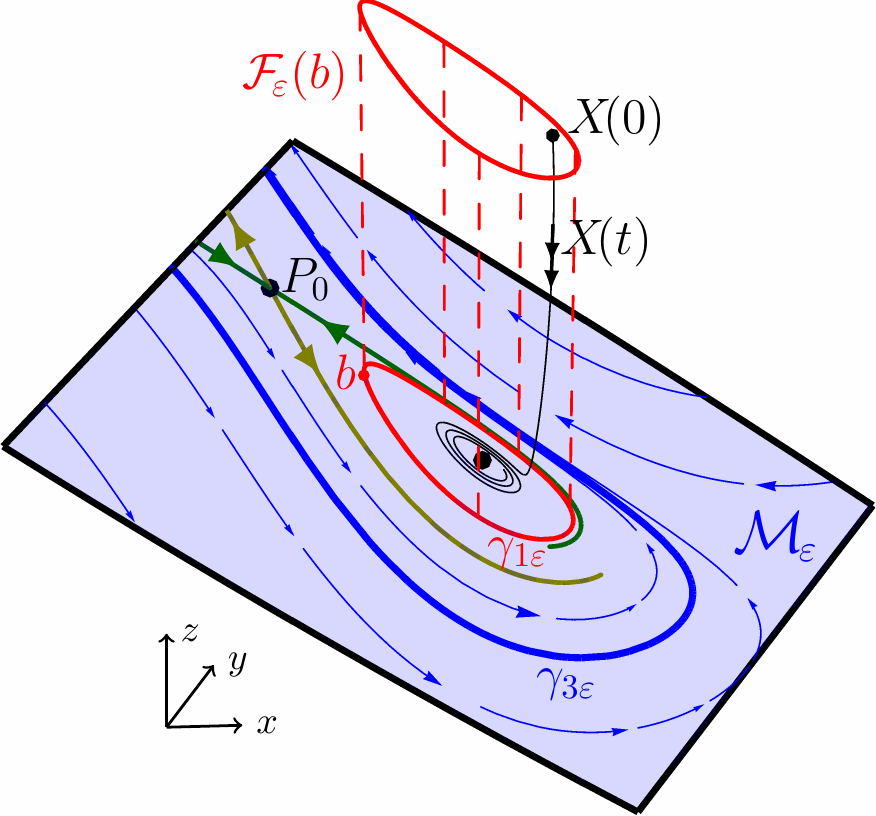}
\caption{
Illustration of the basin of attraction $\mathcal{B}_1$ of the stable 
equilibrium $P_1$
on $\mathcal{M}_\eps$ for $(p,r)$  in the region
between the Hopf bifurcation curve e1-e2 and the homoclinic bifurcation
curve.
Also shown are the unstable limit cycle $\gamma_{1,\eps}$, 
the stable and unstable manifolds of the saddle $P_0$,
the large stable limit cycle $\gamma_{3,\eps}$, 
some fast stable fibers, and a solution $X(t)$ in $\mathcal{B}_1$. 
The union of the fast stable fibers with base points on $\gamma_{1,\eps}$
forms the boundary of this basin.
\label{f-FSsym-BasinRegionIII}
}
\end{figure}

The second case is shown in Figure~\ref{f-FSsym-BasinRegionIII}.
Here, $(p,r)$ is  in the region between
the Hopf bifurcation curve e1-e2 and the homoclinic bifurcation curve; $P_0$ is a saddle, $P_1$ and $P_2$ are stable spirals,
each surrounded by an unstable limit cycle. The limit cycle around $P_1$
is labeled $\gamma_{1,\eps}$; the limit cycle $\gamma_{2,\eps}$ 
around $P_2$ is not shown.
Also, there is a large stable limit cycle, labeled $\gamma_{3, \eps}$, 
which surrounds $P_0, P_1,$ and $P_2$, and the unstable limit cycles.

There are three primary basins of attraction, one for each of the stable
equilibria $P_1$ and $P_2$, and one for the large stable limit cycle $\gamma_{3, \eps}$, denoted by $\mathcal{B}_1$, $\mathcal{B}_2$, and $\mathcal{B}_3$, respectively.
The sets of initial conditions in these primary basins are determined as follows.

First, we observe that all of the initial conditions on the slow manifold
$\mathcal{M}_\eps$ that lie inside the unstable limit cycle $\gamma_{1, \eps}$
are in $\mathcal{B}_1$.  
Denote this set of initial conditions by $\mathcal{B}_1 \vert_{\mathcal{M}_\eps}$.
Then, any initial condition $X(0)$ not on $\mathcal{M}_\eps$ is in $\mathcal{B}_1$
if it lies on the fast stable fiber $\mathcal{F}_\eps(b(0))$ of an initial condition $b(0)$
in $\mathcal{B}_1 \vert_{\mathcal{M}_\eps}$. 
This completely determines the basin~$\mathcal{B}_1$,
see Figure~\ref{f-FSsym-BasinRegionIII}.

Similarly, the basin of attraction $\mathcal{B}_2$ 
of the stable equilibrium $P_2$ consists of 
(i)~the set of all initial conditions on $\mathcal{M}_\eps$
that lie inside the unstable limit cycle $\gamma_{2, \eps}$, 
a set which we denote by $\mathcal{B}_2 \vert_{\mathcal{M}_\eps}$, and
(ii)~the set of all initial conditions $X(0)$ (not on $\mathcal{M}_\eps$) 
that lie on fast stable fibers $\mathcal{F}_\eps(b(0))$ with base points $b(0)$
inside $\mathcal{B}_2 \vert_{\mathcal{M}_\eps}$.

Finally, the basin of attraction $\mathcal{B}_3$ of the large stable limit
cycle $\gamma_{3, \eps}$ consists of 
(i)~the set of all initial conditions on $\mathcal{M}_\eps$ that are exterior
to $\gamma_{1, \eps}$ and $\gamma_{2, \eps}$ and that do not lie
on the stable manifold of $P_0$, a set which we denote by
$\mathcal{B}_3\vert_{\mathcal{M}_\eps}$, and
(ii)~the set of all initial conditions $X(0)$ not on $\mathcal{M}_\eps$ 
which lie on fast stable fibers whose base points are in~$\mathcal{B}_3 \vert_{\mathcal{M}_\eps}$.
 
The second case is more complicated than the first due to the presence of the unstable saddle at $P_0$.
However, the two cases are representative, and the
methods are similar for finding the basin of attraction of any other limit cycle on $\mathcal{M}_{\eps}$.

\section{The Slow--Fast Asymmetric System \label{s-FSasym}}
In this section, we study the effect of asymmetry ($s>0$) in \eqref{MS-xyz} in the limit of large $q$,
\begin{equation}
\begin{split}
\dot{x} &= - x - y , \\
\dot{y} &= ry - pz + s z^2 - y z^2, \\
\eps \dot{z} &= - x - z ,
\end{split}
\Label{FSasym-xyz}
\end{equation}
where $\eps = 1/q$, as before. We refer to \eqref{FSasym-xyz} as the slow--fast asymmetric model. Just as was the case for the symmetric model \eqref{FSsym-xyz}, this slow--fast system has a family of slow invariant manifolds, on which the long--term system dynamics occur and on which the limit cycles lie. By studying the flow on these slow manifolds, we establish the existence of two nondegenerate BT points. For each $s > 0$, the two nondegenerate BT points are organizing centers 
in the $(p,r)$ plane. 
The geometry of these organizing centers and bifurcation curves (and hence also the geometry of the regions where the limit cycles lie) may be understood as a symmetry-breaking of the lone $\mathbb{Z}_2$-symmetric BT point
studied in Section~\ref{ss-FSsymBT}.

\subsection{Slow Manifolds\label{ss-FSasymSlowMan}}
With $\eps=0$, the slow--fast asymmetric model has the same normally hyperbolic critical manifold $\mathcal{M}_0 = \{ (x, y, z): z = - x\}$ as the slow--fast symmetric model~\eqref{FSsym-xyz}. This critical manifold persists for all $0 < \eps \ll 1$. In particular, there is a family of invariant slow manifolds $\mathcal{M}_\eps$, which are given to all orders by the graph of a function $h_\eps(x,y)$, 
\begin{equation}
h_{\eps}(x,y) = h_0(x,y) + \eps h_1(x,y) + \eps^2 h_2(x,y) + \eps^3 h_3(x,y) + \mathcal{O}(\eps^4),
\Label{hepsasymm}
\end{equation}
and this function depends on~$s$. The first two coefficients, $h_0$ and $h_1$, in the expansion of $h_\eps$ are the same as in~(\ref{FSsym-hi}). Then, one finds
\begin{align*}
h_2 &= - ( x + y ) + ( ry + px + (s-y) x^2 ) , \\
h_3 &= - ( 1 - 2p + 4x (y-s) )( x + y ) + (1-r + x^2 )( ry + px + (s-y) x^2 ) .
\end{align*}
On the slow manifolds $\mathcal{M}_{\eps}$, the system \eqref{FSasym-xyz} reduces to the planar system
\begin{equation}		\Label{FSasym-xyheps}
\begin{split}
\dot{x} &= - x - y , \\
\dot{y} &= ry - p h_{\eps}(x,y) + s \left(h_{\eps}(x,y) \right)^2 - y \left( h_{\eps}(x,y) \right)^2.
\end{split}
\end{equation}

To simplify the presentation, we focus on the dynamics of \eqref{FSasym-xyheps} on the critical manifold. That is, we represent $h_\eps$ by the leading order term $h_0=-x$, which gives 
\begin{equation}		\Label{FSasym-xy}
\begin{split}
\dot{x} &= - x - y , \\
\dot{y} &= ry + px + s x^2 - y x^2.
\end{split}
\end{equation}
The dynamics of \eqref{FSasym-xyheps} for $0<\eps \ll 1$ are regular perturbations of those of \eqref{FSasym-xy}.

\subsection{Organizing Centers\label{ss-FSasymBT}}
\begin{figure}[h!]
\centering
\includegraphics[width=3.5in]{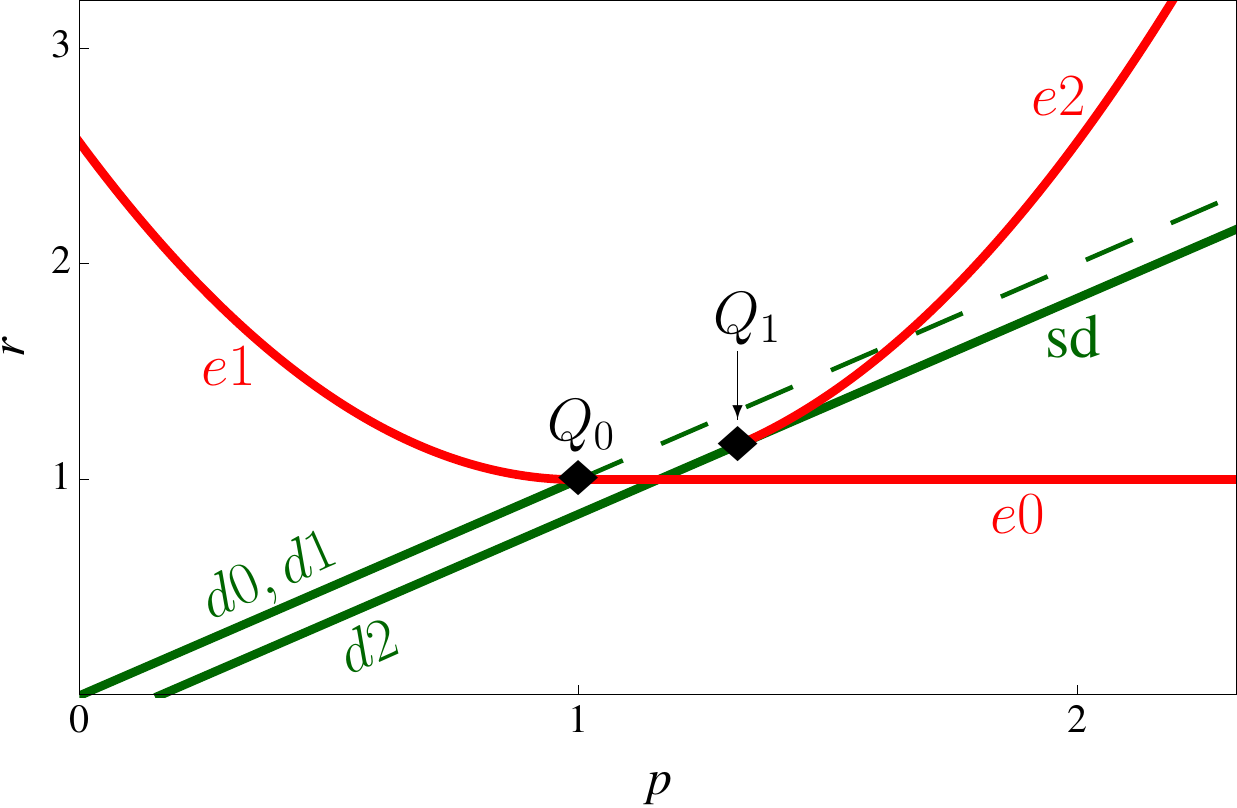}
\caption{Stability regions of $P_0$, $P_1$, and $P_2$ for \eqref{FSasym-xy} with $s = 0.8$.
$P_0$ is linearly stable in the region enclosed by d0, e0, and the $p$ axis. On e0, $P_0$ undergoes a supercritical Hopf bifurcation
with natural frequency $\omega_0^* = \sqrt{p-1}$.
$P_1$ is linearly stable in the region enclosed by d1, e1, and the $r$ axis.
On e1, $P_1$ loses stability due to a subcritical Hopf bifurcation with natural frequency $\omega_1^* = \sqrt{2(r-1) + s \sqrt{r-1}}$.
$P_2$ is linearly stable in the region enclosed by the concatenation of d2 and e2 and the $r$ axis.
On e2, $P_2$ loses stability due to a subcritical Hopf bifurcation with natural frequency $\omega_2^* =  \sqrt{2(r-1) - s \sqrt{r-1}}$.
The diagonal is the curve of transcritical bifurcations of $P_0$ and $P_1$. The shifted diagonal is the curve of saddle-node bifurcations, where $P_1$ and $P_2$ coalesce.
\label{f-2Dasym-stability}}
\end{figure}

The origin $P_0 = (0,0)$ is an equilibrium of \eqref{FSasym-xy} for all $p$, $r$, and $s$. If $r > p - \tfrac14 s^2$, there are two additional equilibria, namely $P_1 =  (x_1^*, - x_1^*)$ and $P_2 =  (x_2^*, - x_2^*)$, where
\begin{equation}
x_1^* = \tfrac12 [- s + \sqrt{s^2 + 4(r - p)}] , \quad
x_2^* = \tfrac12 [- s - \sqrt{s^2 + 4(r - p)}] .
\Label{FSasym-x1*x2*}
\end{equation}
Figure \ref{f-2Dasym-stability} shows the results of a linear stability analysis of the equilibria of \eqref{FSasym-xy}.

We refer to the line $r = p - \tfrac14 s^2$ as the shifted diagonal (marked d2 and sd in Figure~\ref{f-2Dasym-stability}). 
Note that $x_2^* < x_1^* < 0$ if $p - \tfrac14 s^2 < r < p$, and $x_2^* < 0 < x_1^*$ if $r > p$.

At the points
\begin{equation} 
Q_0 = (1, 1) \quad \text{ and } \quad Q_1 = (1 + \tfrac12 s^2, 1 + \tfrac14 s^2),
\Label{Q1Q2} 
\end{equation}
the system \eqref{FSasym-xy} has a zero eigenvalue of geometric multiplicity two (at $P_0$ and $P_1$, respectively). Thus, $Q_0$ and $Q_1$ are organizing centers.

\begin{figure}[h!]
\centering
\includegraphics[width=5in]{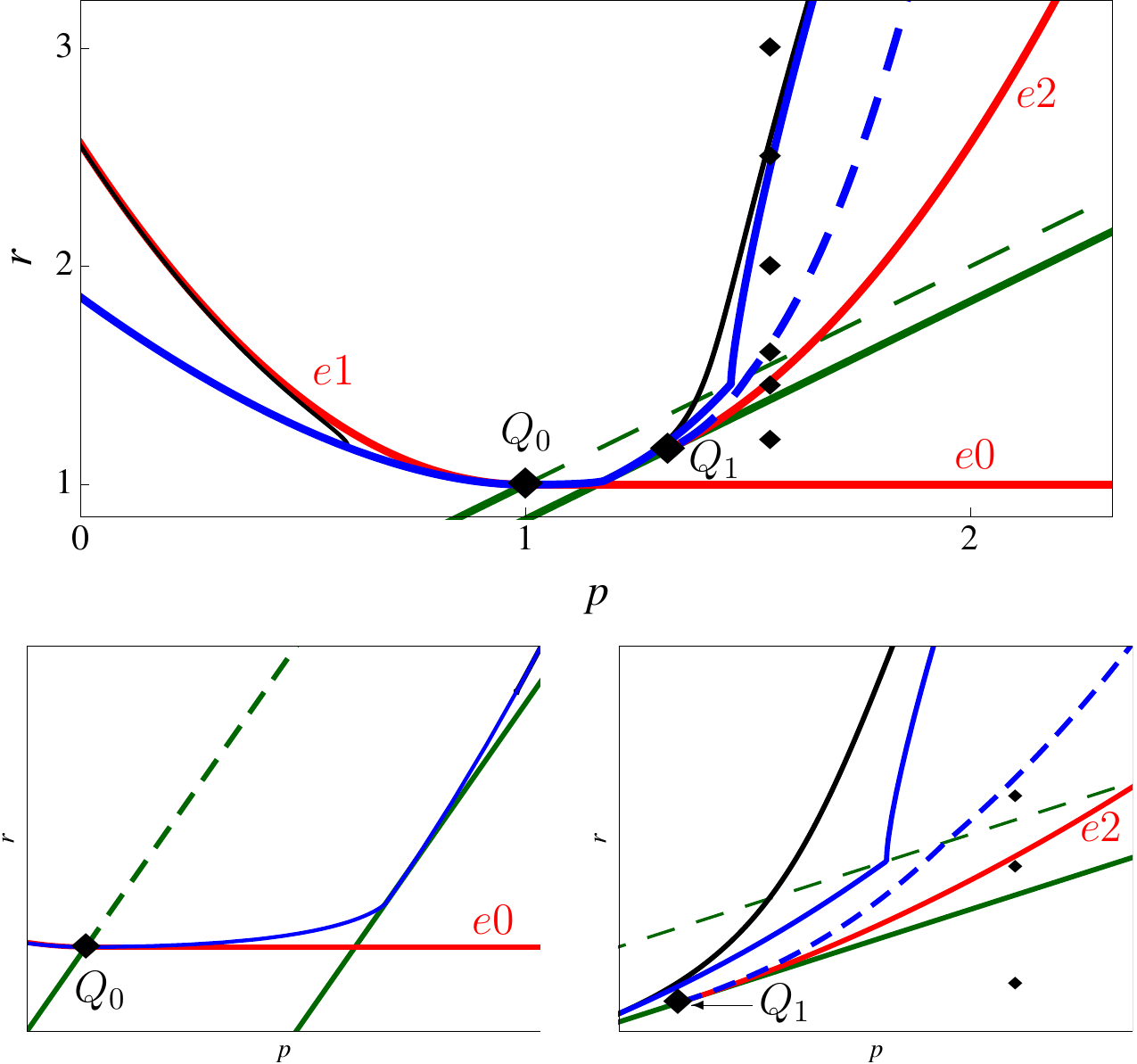}
\put(-362,327){(a)}
\put(-370,124){(b)}
\put(-182,124){(c)}
\caption{(a) Local and global bifurcation curves of system~\eqref{FSasym-xy} with $s = 0.8$. 
The labels are as in Figure \ref{f-2Dasym-stability}.
Also shown are magnified views of the neighbourhoods of (b) $Q_0$ and (c) $Q_1$. 
\label{f-2Dasym-BifurcationCurves}}
\end{figure}

\begin{figure}[h!]
\centering
\includegraphics[width=5in]{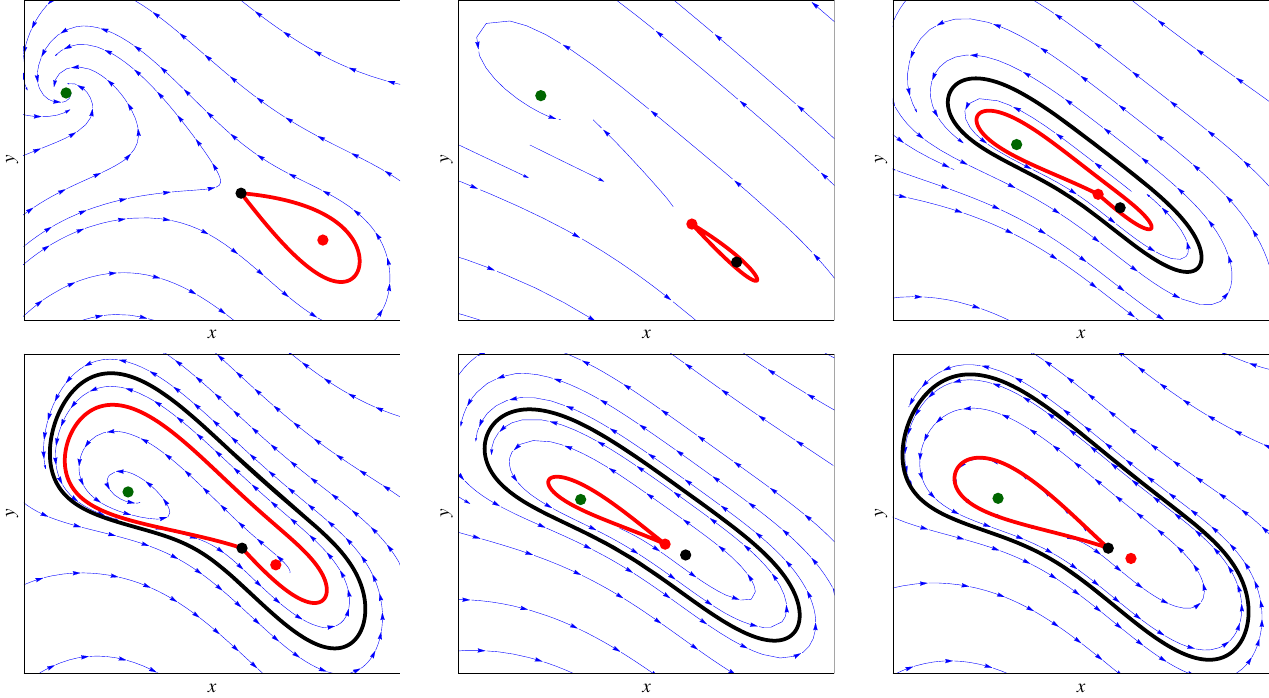}
\put(-371,190){(a)}
\put(-248,190){(b)}
\put(-123,190){(c)}
\put(-371,88){(d)}
\put(-247,88){(e)}
\put(-123,88){(f)}
\caption{Homoclinic orbits (red curves) of \eqref{FSasym-xy}. The black, red and green markers correspond to $P_0, P_1$ and $P_2$, respectively. 
(a)~Right homoclinic to $P_0$ for $(p,r) \approx (0.35, 1.40631)$.
(b)~Right homoclinic to $P_1$ for $(p,r) \approx (1.1,1.00308)$.
(c)~Large-amplitude homoclinic to $P_1$ for $(p,r) \approx (1.37,1.27743)$. 
(d)~Large-amplitude homoclinic to $P_0$ for $(p,r) \approx (1.49, 1.87285)$.
(e)~Left homoclinic to $P_1$ for $(p,r) \approx (1.45,1.36238)$. 
(f)~Left homoclinic to $P_0$ for $(p,r) \approx (1.6, 1.80921)$.
In frames (c)--(f), the black curve is the large-amplitude stable limit cycle.
\label{f-2Dasym-Homoclinics}}
\end{figure}
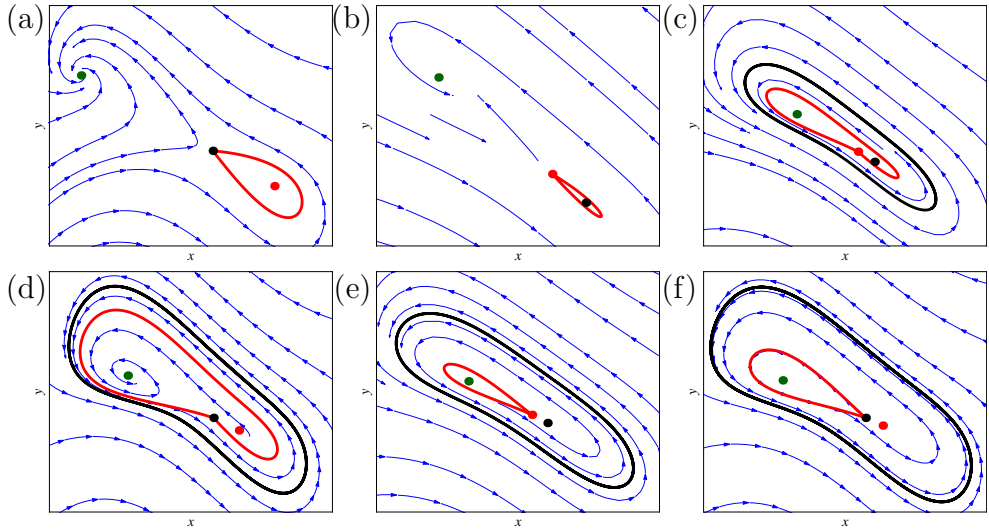

Figure \ref{f-2Dasym-BifurcationCurves} shows the organizing centers and the branches of global bifurcations of \eqref{FSasym-xy}.
There are two branches of homoclinic bifurcations emanating from $Q_0$ (solid blue), one to the left and one to the right, and a single branch emanating from $Q_1$ to the right (dashed blue). 
The companion Figure \ref{f-2Dasym-Homoclinics} shows the types of homoclinic orbits along these branches. 
In Figure~\ref{f-2Dasym-Homoclinics}, the equilibria $P_0, P_1,$ and $P_2$ are marked by black, red, and green dots, respectively.  
These branches, and the homoclinic orbits along them, are as follows:

\begin{itemize}
\setlength{\itemsep}{0em}
\item The branch emanating from $Q_0$ to the left consists of right-homoclinic orbits to $P_0$, which enclose $P_1$. A representative orbit is shown in Figure \ref{f-2Dasym-Homoclinics}(a).
\item The branch emanating from $Q_0$ to the right consists of three segments:
\begin{itemize}
\item A segment from $Q_0$ to the shifted diagonal, most clearly seen in Figure \ref{f-2Dasym-BifurcationCurves}(b); this branch consists of right-homoclinic orbits to $P_1$, which enclose $P_0$. A representative orbit is shown in Figure \ref{f-2Dasym-Homoclinics}(b). 
\item A segment from the shifted diagonal to the diagonal, most clearly seen in Figure \ref{f-2Dasym-BifurcationCurves}(a) and (c); this branch consists of large-amplitude double-loop homoclinic orbits to $P_1$, which enclose $P_0$ and $P_2$. A representative orbit is shown in Figure \ref{f-2Dasym-Homoclinics}(c). 
\item A segment beyond the diagonal, most clearly seen in Figure \ref{f-2Dasym-BifurcationCurves}(a); this branch consists of large-amplitude double-loop homoclinic orbits to $P_0$, which enclose $P_1$ and $P_2$. A representative orbit is shown in Figure \ref{f-2Dasym-Homoclinics}(d). 
\end{itemize}
\item The branch of homoclinic bifurcations emanating from $Q_1$ (dashed blue) is made up of two segments:
\begin{itemize}
\item A segment from $Q_1$ to the diagonal, most clearly seen in Figure \ref{f-2Dasym-BifurcationCurves}(c); this branch consists of left-homoclinics to $P_1$, which enclose $P_2$. A representative orbit is shown in Figure \ref{f-2Dasym-Homoclinics}(e).
\item A segment beyond the diagonal, most clearly seen in Figure \ref{f-2Dasym-BifurcationCurves}(a); this branch consists of left-homoclinics to $P_0$, which enclose $P_2$. A representative orbit is shown in Figure \ref{f-2Dasym-Homoclinics}(f).
\end{itemize}
\end{itemize}
The existence of all six types of homoclinics will be formally proven in Section \ref{ss-BTanalysis}.
Figure \ref{f-2Dasym-BifurcationCurves} also shows two branches of saddle-node bifurcations of limit cycles (black curves). We note that these curves are no longer attached to the organizing centers due to the asymmetry, i.e., since the BT points are non-degenerate.
The local and global bifurcation curves serve as boundaries between the regions of different dynamical behaviour. 

\begin{figure}[h!]
\centering
\includegraphics[width=5in]{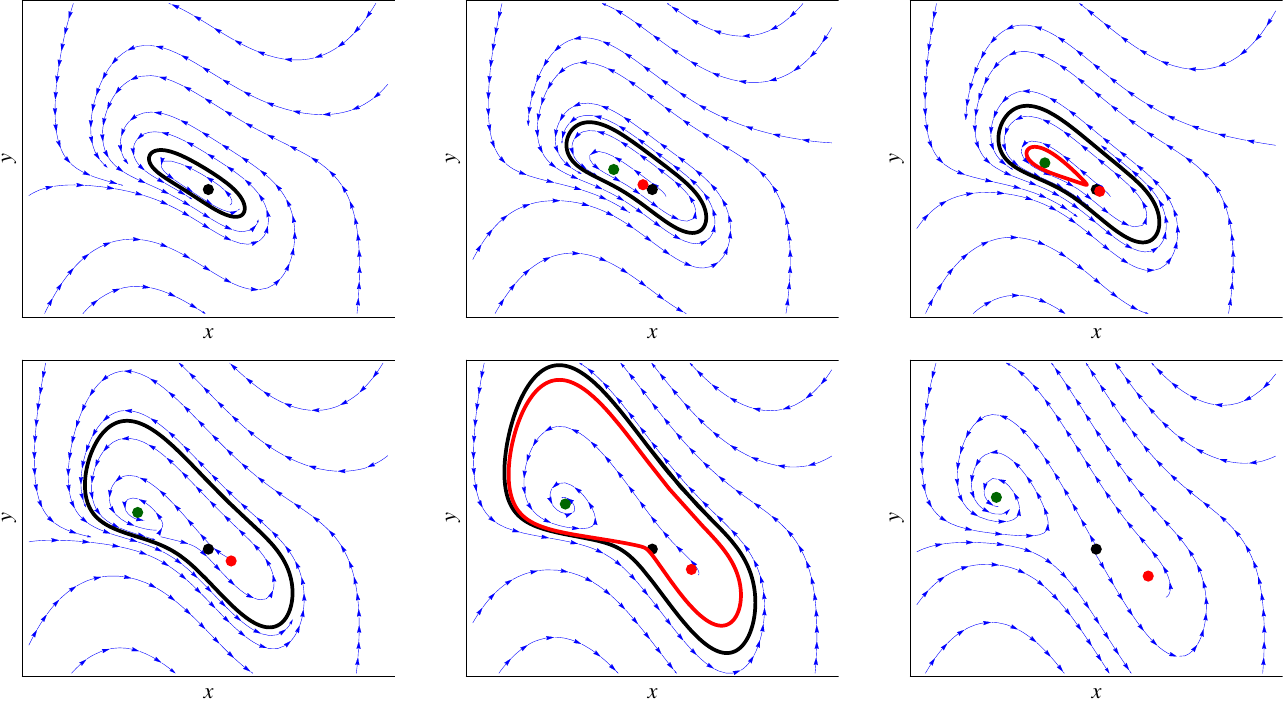}
\put(-372,188){(a)}
\put(-248,188){(b)}
\put(-121,188){(c)}
\put(-372,86){(d)}
\put(-248,86){(e)}
\put(-121,86){(f)}
\caption{Phase planes corresponding to the small, black diamond markers
in Figure~\ref{f-2Dasym-BifurcationCurves} along the vertical line $p=1.55$.
(a)~For $r=1.2$, there is a stable limit cycle (black curve)  surrounding $P_0$ (black dot).
(b)~For $r=1.45$, the equilibria $P_1$ and $P_2$ (red and green dots) exist.
(c)~For $r=1.6$, there is an unstable limit cycle (red) surrounding $P_2$, and
the relative positions of $P_1$ and $P_0$ have switched. 
(d)~For $r=2$, the unstable limit cycle has disappeared in the homoclinic bifurcation.
(e)~For $r=2.5$, a large-amplitude unstable limit cycle (red curve) exists
inside the large-amplitude stable limit cycle.
(f)~For $r=3$, the equilibrium $P_2$ is the only attractor, since the large-amplitude stable
and unstable limit cycles have disappeared in a saddle-node bifurcation.
\label{f-2Dasym-PhasePlanes}}
\end{figure}
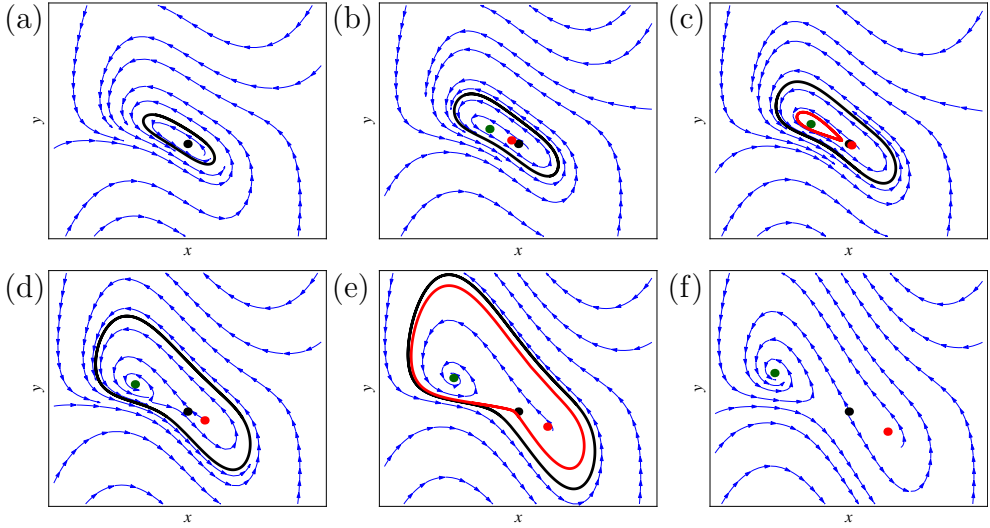

Figure~\ref{f-2Dasym-PhasePlanes} shows the sequence of phase portraits
that are encountered for a fixed value of $p>1+\tfrac12 s^2$
as $(p, r)$ moves through the different regions for increasing values of $r$. 
There is one phase portrait for each of the six small black diamonds
in Figure~\ref{f-2Dasym-BifurcationCurves}. 

In Figure~\ref{f-2Dasym-PhasePlanes}(a), the parameters are chosen in the region bounded by e0 and the shifted diagonal, where $P_0$ is an unstable spiral surrounded by a stable limit cycle (black curve). In frame~(b), $P_1$ (red dot) and $P_2$ (green dot) lie inside the large stable limit cycle, for $(p,r)$ values in the region between the shifted diagonal and e2. In frame~(c), we see the unstable limit cycle (red) surrounding $P_2$ that lies in the region bounded by e2 and the dashed curve of homoclinics. Next, frame~(d) shows a representative phase portrait obtained after the homoclinic bifurcation curve (dashed blue curve in Figure~\ref{f-2Dasym-BifurcationCurves}) is crossed. One sees that all solutions (not on the stable manifold of $P_0$) are forward asymptotic to $P_2$ (green dot) or to the large stable limit cycle (black curve). In frame~(e), there are large-amplitude unstable (red) and stable (black) limit cycles, in the narrow region between the homoclinic bifurcation curve (blue) and the curve of saddle-node bifuractions of limit cycles (black).
These disappear as one crosses the curve of saddle-node bifurcations of limit cycles (upper black curve in Figure~\ref{f-2Dasym-BifurcationCurves}). In frame~(f), only the three equilibria remain.

\begin{remark}
The bifurcation structure of \eqref{FSasym-xy} collapses to that of the symmetric case \eqref{FSsym-xyheps} as $s \to 0$ (see Figure \ref{f-2Dasym-SymmetryBreaking}).  
More specifically, the shifted diagonal line of saddle-node bifurcations of $P_1$ and $P_2$ collapses onto the diagonal line of transcritical bifurcations of $P_0$ and $P_1$, thus creating the diagonal line of pitchfork bifurcations of system \eqref{FSsym-xyheps}.
The organizing centers $Q_0$ and $Q_1$ merge and become the $\mathbb{Z}_2$-symmetric BT point at $Q = (1,1)$ as $s \to 0$.
Concomitantly, the curves, e1 and e2, of Hopf bifurcations merge to the single curve labelled e1-e2 in Figure \ref{f-FSsym-HeatPlot}. 
Additionally, the curves of homoclinic bifurcations merge into a single curve (which is tangent to the line $r-1= -4(p-1)$ at $Q$) as $s \to 0$. 
Similarly, the curves of saddle-node bifurcations of limit cycles collapse to a single curve (which is tangent to the line $r-1 \approx -3.03(p-1)$ at $Q$) as $s \to 0$.  

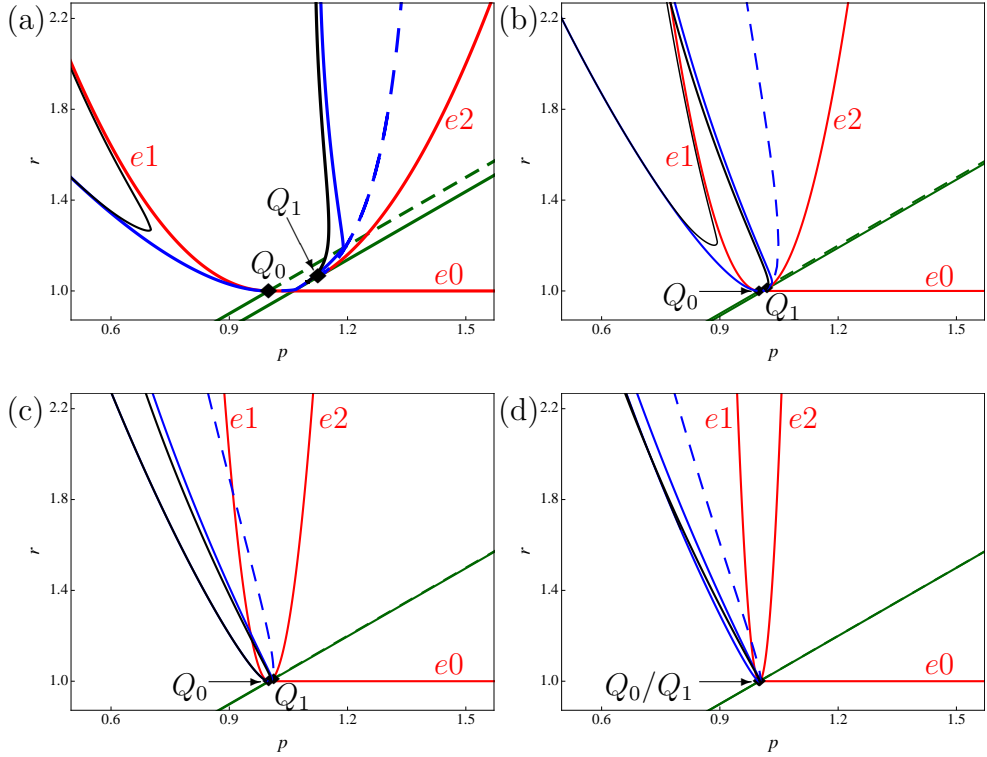
\begin{figure}[h!]
\centering
\includegraphics[width=5in]{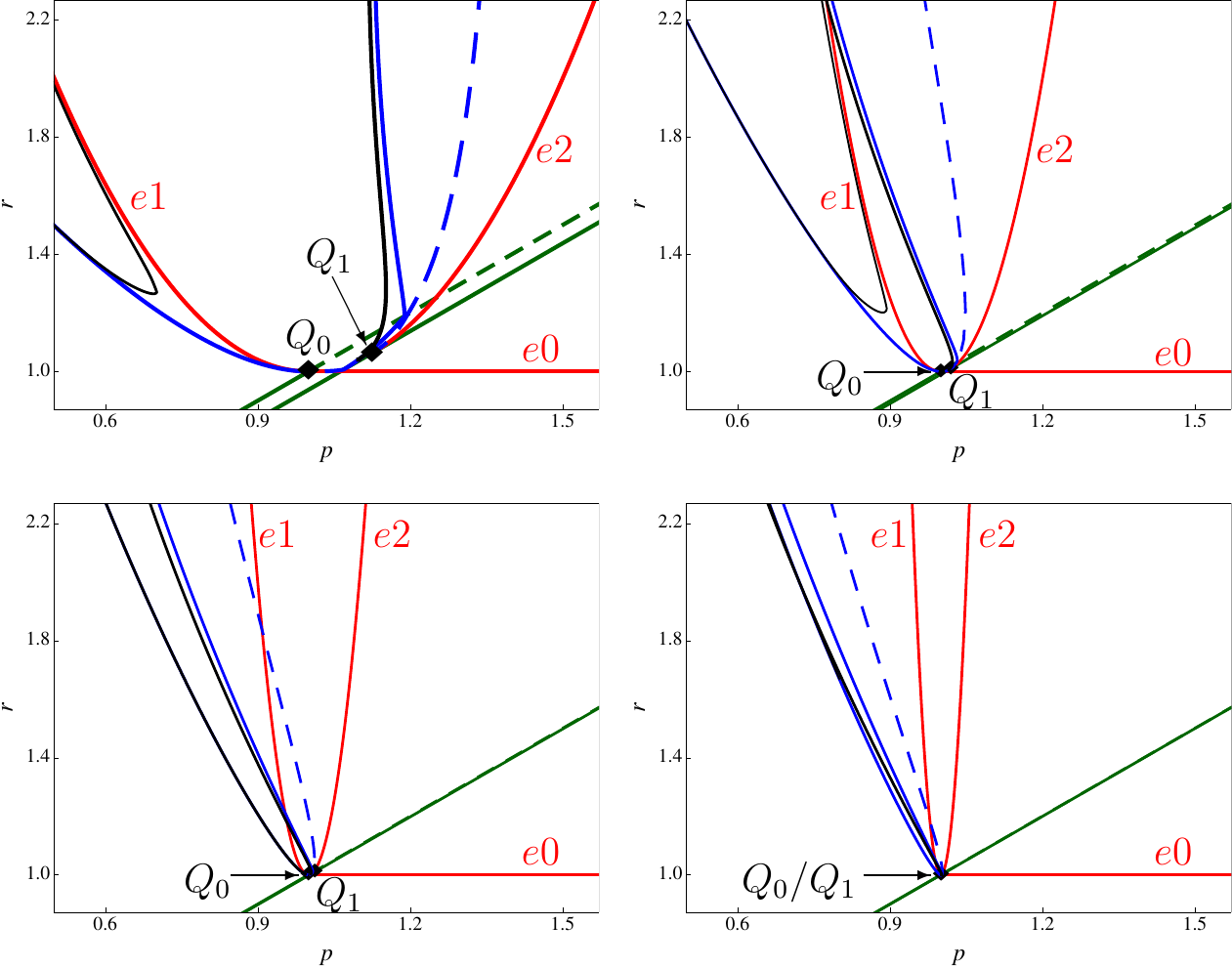}
\put(-370,274){(a)}
\put(-184,274){(b)}
\put(-370,126){(c)}
\put(-184,126){(d)}
\caption{Bifurcation structure of system~\eqref{FSasym-xy} for (a) $s = 0.5$, (b) $s=0.2$, (c) $s=0.1$ and (d) $s = 0.05$.
This sequence illustrates how the bifurcation curves for $s>0$ collapse onto those of \eqref{FSsym-xyheps} as $s \to 0$, recall Figure \ref{f-FSsym-HeatPlot}. The colors are the same as in Figure \ref{f-2Dasym-BifurcationCurves}.
\label{f-2Dasym-SymmetryBreaking}}
\end{figure}

\end{remark}

\subsection{Bogdanov-Takens Unfolding Analysis} \label{ss-BTanalysis}
\begin{figure}[h!]
\centering
\includegraphics[width=5in]{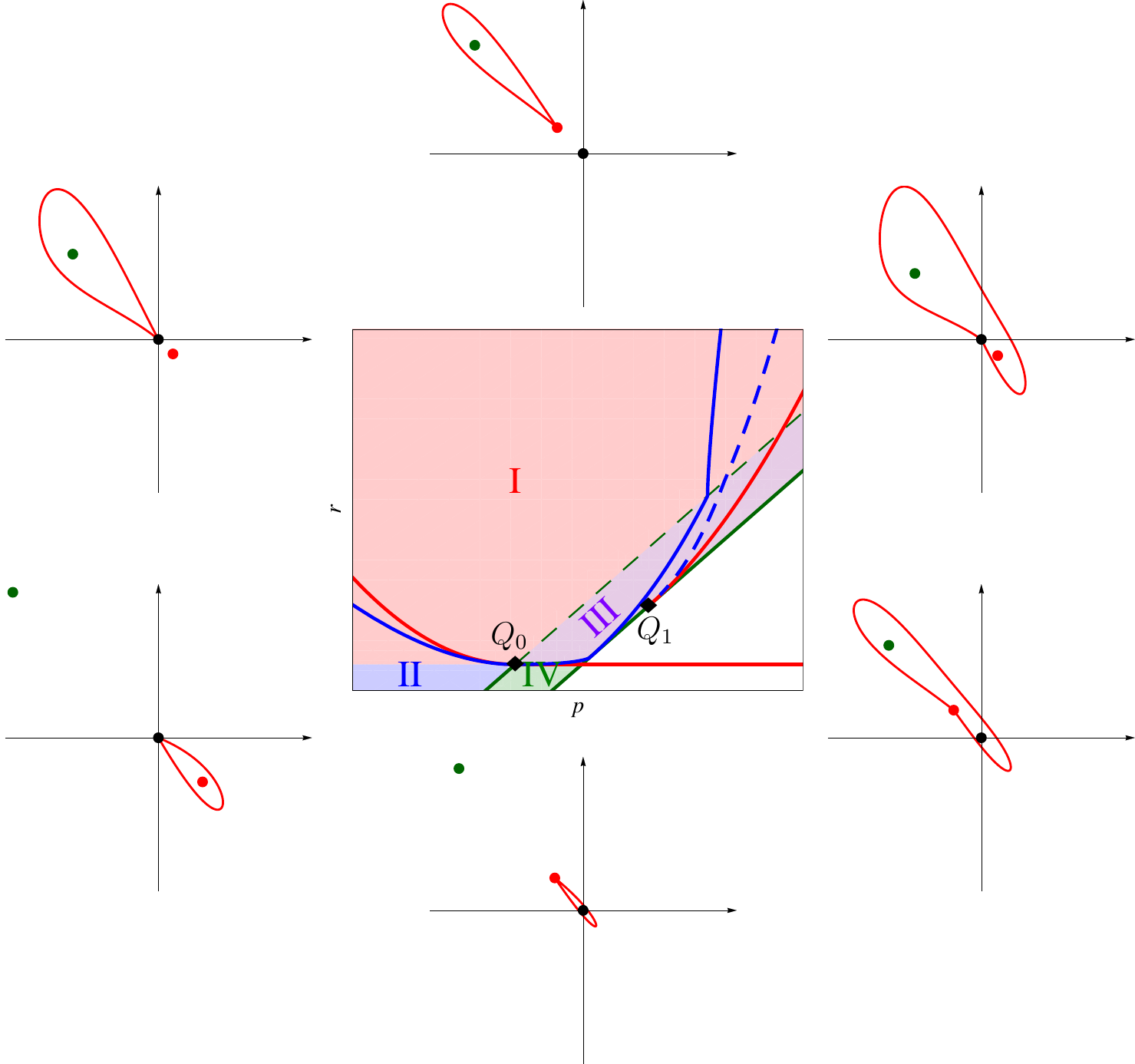}
\caption{Regions I--IV in the $(p,r)$ plane (central frame) and the different types of homoclinic orbits (red curves) in the original $(x,y)$ phase plane (surrounding frames). The black, red, and green dots correspond to the equilibria $P_0, P_1,$ and $P_2$, respectively. Note that the curves of saddle-node bifurcations of limit cycles have been omitted from the $(p,r)$ plane.}
\label{f-2Dasym-regions}
\end{figure}

In this section, we present the unfolding analysis of the non-degenerate BT points $Q_0$ and $Q_1$ in system \eqref{FSasym-xy}. Specifically, we study the homoclinic orbits of an appropriate Hamiltonian system and use Melnikov theory \cite{Guckenheimer1983,Melnikov1963} to determine the parameter sets for which these homoclinics persist under small perturbations. In this manner, we formally prove the persistence of the six types of homoclinic orbits, which lie along the blue branches that emanate from $Q_0$ and $Q_1$ in Figure \ref{f-2Dasym-Homoclinics}. Figure \ref{f-2Dasym-regions} summarizes the results of this subsection.

\subsubsection{Rescaling and Partition of the Parameter Plane\label{sss-rescaling}}
First, we make the change of variables $(x, -(x+y)) \mapsto (x, y)$, so that \eqref{FSasym-xy} becomes
\begin{equation}
\begin{split}
\dot{x} &= y , \\
\dot{y} &= (r - p) x + (r - 1) y - (s + y) x^2 - x^3 .
\end{split}
\Label{FSasym-xy'}
\end{equation}
The equilibria are $P_0 = (0, 0)$, $P_1 = (x_1^*, 0)$, and $P_2 = (x_2^*, 0)$, where $x_1^*$ and $x_2^*$ are again given by~(\ref{FSasym-x1*x2*}).
We rescale the dependent and independent variables and the parameters by 
\begin{equation}
x(t) = \eta\, u({\tilde t}) , \quad
y(t) = \eta^2\, v ({\tilde t}), \quad
{\tilde t} = \eta\, t, \quad
\mu = \frac{r - p}{\eta^2} , \quad
\lambda = \frac{r - 1}{\eta^2} , \quad
\delta = \frac{s}{\eta}.
\Label{FSasym-rescaling}
\end{equation}
With these rescaled variables and parameters, system \eqref{FSasym-xy'} is equivalent with
\begin{equation}
\begin{split}
\dot{u} &= v , \\
\dot{v} &= \mu u - \delta u^2 - u^3 + \eta (\lambda - u^2) v ,
\end{split}
\Label{FSasym-uv}
\end{equation}
where the overdot denotes $\frac{d}{d{\tilde t}}$, and we drop the tildes. 

There are four distinct regions in the $(p,r)$ plane, depending on $\mu$ and $\lambda$ (Figure \ref{f-2Dasym-regions}):
\begin{enumerate}[I.]
\setlength{\itemsep}{0pt}
\item The set $\{ \mu>0, \, \lambda>0 \}$ corresponds to the region of the $(p,r)$ plane above the main diagonal and above the line $r=1$ (red shaded region).
\item The set $\{ \mu>0, \,\, \lambda<0 \}$ is the region above the main diagonal and below the line $r=1$ (blue shaded region).
\item The set $\{ -\tfrac{1}{4} \delta^2 < \mu<0 , \,\, \lambda>0 \}$ is the region enclosed by the two diagonals, and lies above the line $r=1$ (purple shaded region). 
\item The set $\{ -\tfrac{1}{4} \delta^2 <\mu<0, \,\, \lambda<0 \}$ corresponds to the region enclosed by the two diagonals and lies below the line $r=1$ (green shaded region).
\end{enumerate}
The organizing centers are $Q_0=(0,0)$ and $Q_1=(-\tfrac{1}{4}\delta^2,\tfrac{1}{4}\delta^2)$ in the $(\mu,\lambda)$ plane, corresponding to $P_0=(0,0)$ and $P_1=(-\tfrac{1}{2}\delta,0)$, respectively.

\subsubsection{Persistence of the Homoclinics to $P_0$\label{sss-Q1homoclinics}}
To study the homoclinics to $P_0$, we consider the region of $\mu >0$. 
Without loss of generality, we set $\mu = 1$ in \eqref{FSasym-uv} and begin with the $\eta = 0$ limit of \eqref{FSasym-uv}
\begin{equation}
\begin{split}
\dot{u} &= v , \\
\dot{v} &= u - \delta u^2 - u^3,
\end{split}
\Label{FSasym-HamiltonianDS-uv}
\end{equation}
which is Hamiltonian. The Hamiltonian function is $H (u, v) = \tfrac12 v^2 - \tfrac12 u^2 + \tfrac13 \delta u^3 + \tfrac14 u^4 $. System~\eqref{FSasym-HamiltonianDS-uv} has equilibria at $(0,0)$, $(u_1^*, 0)$, and $(u_2^*, 0)$, where
\[
u_1^* (\delta) = \frac{1}{2} \left(- \delta + \sqrt{\delta^2+4} \right) , \quad
u_2^* (\delta) = \frac{1}{2} \left(- \delta - \sqrt{\delta^2+4} \right) .
\]

\begin{figure}[h!]
\centering
\includegraphics[width=3in]{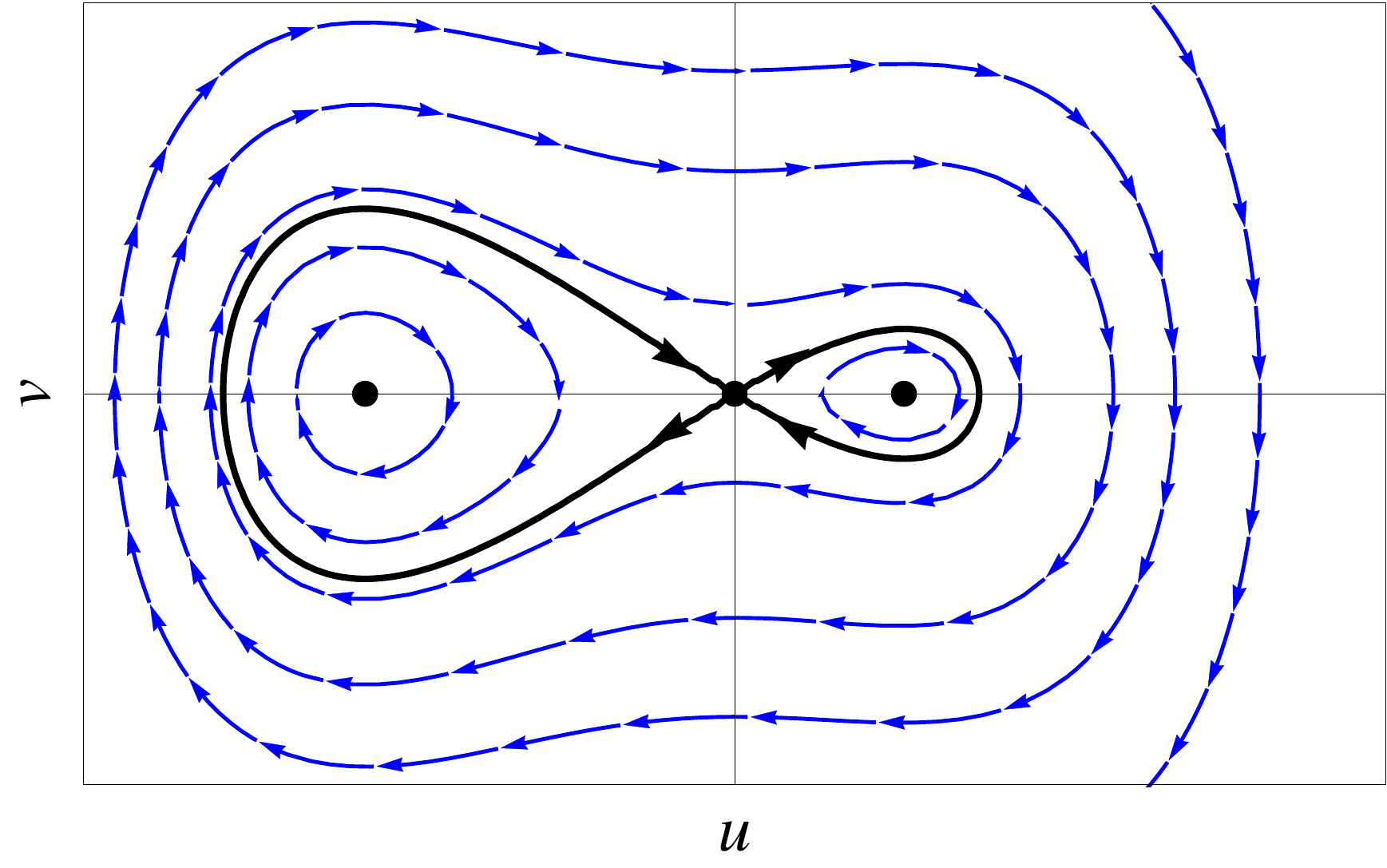}
\caption{Phase portrait of the Hamiltonian system \eqref{FSasym-HamiltonianDS-uv}.
 \label{f-FSasym-PhasePortrait1}}
\end{figure}

The saddle $P_0$ is connected to itself by a pair of asymmetric homoclinic orbits, $\Gamma^-_0$ and $\Gamma^+_0$, which surround $(u_2^*, 0)$ and $(u_1^*, 0)$, respectively. They are given explicitly by 
\begin{equation}
\begin{split}
\Gamma_0^\pm : \quad 
\left( u_0^\pm(t), v_0^\pm(t) \right) &= \left( \frac{\pm 3\alpha}{ {\cosh} (t) \pm \alpha \delta}, \frac{\mp 3\alpha \,{\sinh} (t)} {\left( {\cosh} (t) \pm \alpha \delta \right)^2} \right),
\end{split}
\Label{FSasym-leftrighthomoclinics}
\end{equation}
where $\alpha = \sqrt{ \frac{2}{9 + 2 \delta^2} }$.
The phase portrait of \eqref{FSasym-HamiltonianDS-uv} is shown in Figure~\ref{f-FSasym-PhasePortrait1}, and we note that the solutions $\Gamma_0^\pm$ limit exactly, as $\delta \to 0$, onto the homoclinic orbits of the unperturbed reduced system on $\mathcal{M}_0$ in the slow--fast symmetric system. 

\begin{figure}[h!]
\centering
\includegraphics[width=3.5in]{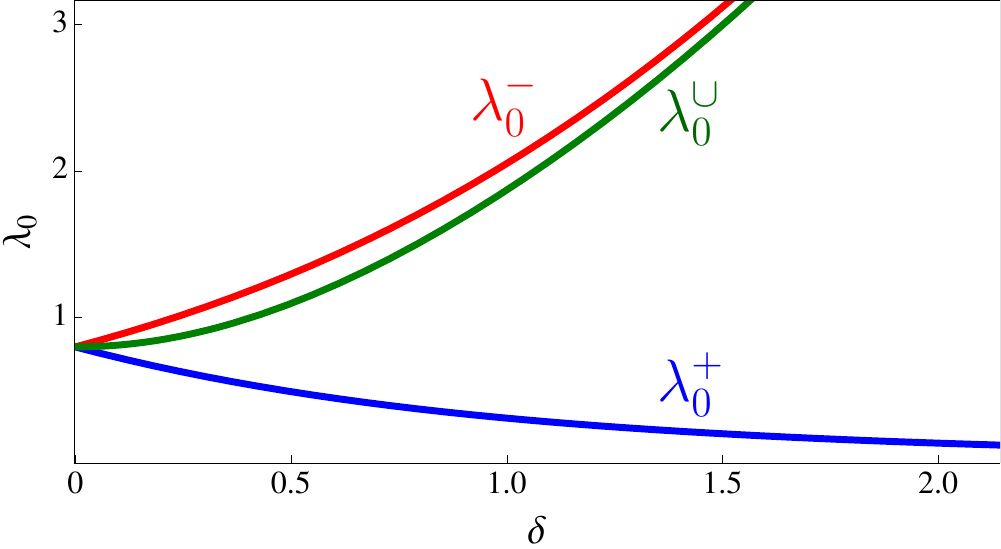}
\caption{Simple zeroes $\lambda_0^{\pm}$ and $\lambda_0^\cup$ of the distance function for $\mu = 1$, corresponding to the right, left, and large-amplitude homoclinics to $P_0$ (see \eqref{eq:lambdapmP0} and \eqref{eq:lambdaLAHP0}). The three branches coalesce to $\lambda=\frac{4}{5}$ as $\delta \to 0$, indicating that the right, left, and large-amplitude homoclinic orbits to $P_0$ collapse onto the simple pair of symmetric homoclinics studied in Section \ref{s-FSsym}.}
\label{f-2Dasym-lambda-p0}
\end{figure}

Now, for each of the unperturbed homoclinic orbits, $\Gamma^+_0$ and $\Gamma^-_0$, we develop the distance function, $D(\lambda)$, as an asymptotic expansion in $\eta$:
\begin{equation}
\begin{split}
D(\lambda) &= \eta \int_{-\infty}^\infty \left( \lambda - (u_0^\pm (t))^2 \right) (v_0^\pm (t))^2 \, dt + \mathcal{O}(\eta^2) 
= \eta \left( \lambda I_0^\pm (\delta) - I_2^\pm (\delta) \right) + \mathcal{O}(\eta^2) ,
\end{split}
\Label{FSasym-Melnikov1}
\end{equation}
where we define the integrals $I_0^\pm(\delta)$ and $I_2^\pm(\delta)$ as
\begin{equation}	
\begin{split}
I_0^\pm (\delta) &= \int_{-\infty}^\infty (v_0^\pm (t))^2 \, dt = \frac{2}{3} \left( 3\phi^2-1 \right)+4 \psi \phi^2 \tan^{-1} \left( \psi \mp \phi \right), \\ 
I_2^\pm (\delta) &= \int_{-\infty}^\infty (u_0^\pm (t) \, v_0^\pm (t))^2 \, dt = \frac{16}{15} + \frac{23}{3} \psi^2 + 7\psi^4+2\psi \phi^2 \left( 3\phi^2+4\psi^2 \right) \tan^{-1} \left( \psi \mp \phi \right),
\end{split}
\Label{eq:MelnikovIntegrals}
\end{equation}
where $\psi = \tfrac{1}{3} \delta \sqrt{2}$ and $\phi = \tfrac{1}{3}\sqrt{2\delta^2+9}$ (and $\phi^2-\psi^2 =1$).
We note that $I_0^\pm(\delta)$ is the area enclosed by $\Gamma_0^\pm$ in the $(u,v)$ plane, since $\int v^2\, dt = \int v \, du$.
The bifurcation equation $D(\lambda) = 0$ gives the parameter set for which the homoclinics persist.
Thus, the homoclinics, $\Gamma_0^\pm$, to $P_0$ persist for
\begin{equation} 
\lambda = \lambda_0^\pm(\delta) = \frac{I_2^\pm (\delta)}{ I_0^\pm (\delta)} + \mathcal{O}(\eta),
\Label{eq:lambdapmP0}
\end{equation}
which is shown in Figure \ref{f-2Dasym-lambda-p0}.  This demonstrates the existence of the curves of homoclinic bifurcations to $P_0$. More specifically, $\lambda_0^+$ corresponds to the solid branch of right-homoclinics in region I that emanates to the left of $Q_0$, and these are the right-homoclinic orbits shown in Figure \ref{f-2Dasym-Homoclinics}(a).
Also, $\lambda_0^-$ corresponds to the dashed branch of left-homoclinics in region I that emanates to the right of $Q_1$ and lies above the main diagonal, and these are the left-homoclinic orbits shown in Figure \ref{f-2Dasym-Homoclinics}(f). 

\begin{figure}[h!]
\centering
\includegraphics[width=3.5in]{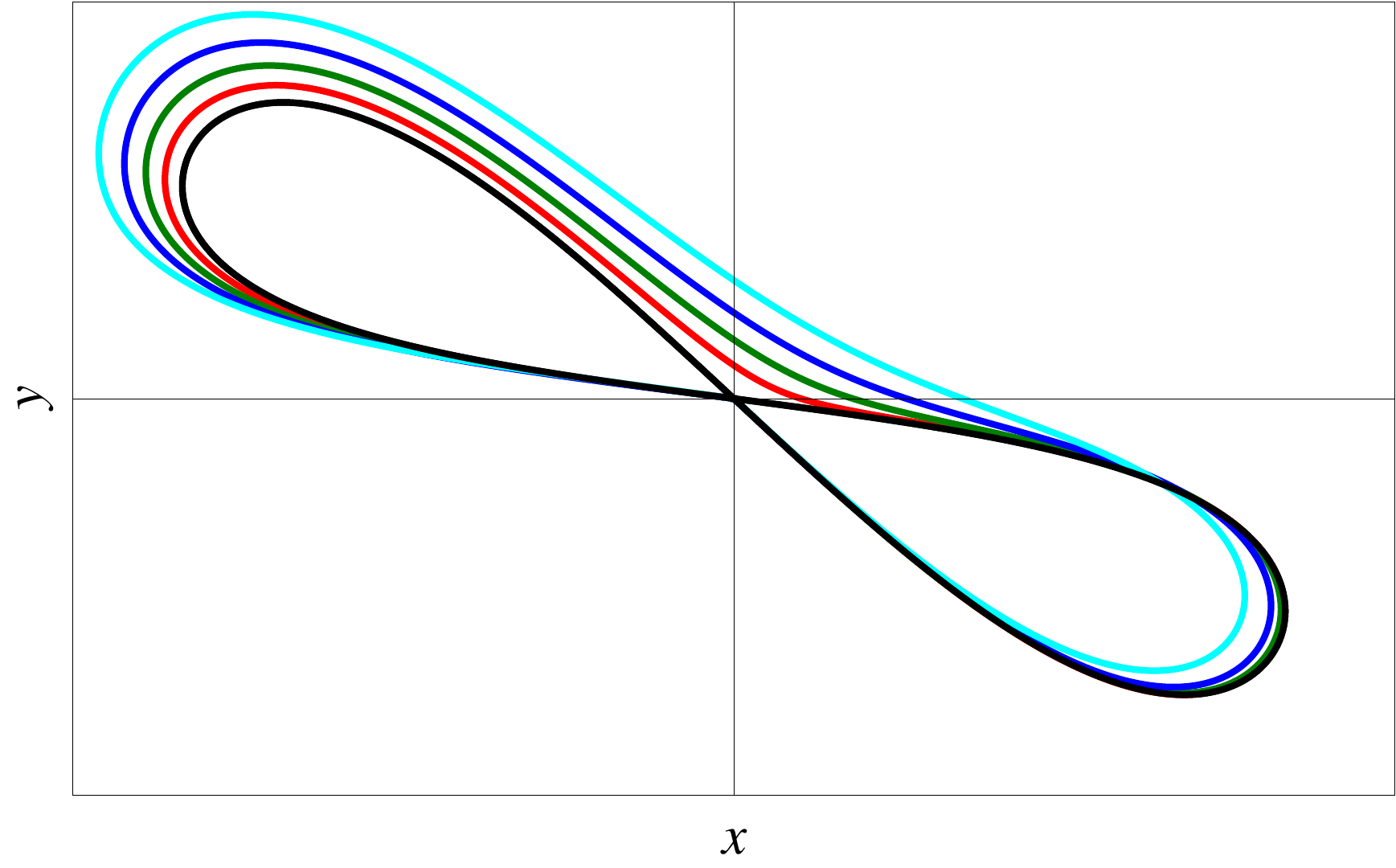}
\caption{Large-amplitude homoclinics to $P_0$ for $s=0.02, 0.05, 0.1$ and $s=0.2$ in red, green, blue, and cyan, respectively. The symmetric limit is shown in black. Note that the values of $p$ and $r$ depend on $s$, and are determined from the numerical continuation.}
\label{f-2Dasym-LAH}
\end{figure}

To establish the existence and persistence of the large-amplitude double-loop homoclinics to $P_0$ (see Figure \ref{f-2Dasym-Homoclinics}(d)), we make the following numerical and analytical observations. 
Numerically, the large-amplitude homoclinics to $P_0$, as illustrated in Figure \ref{f-2Dasym-regions} for $s=0.8$ for instance, converge in the symmetric (i.e., $s \to 0$) limit to the pair of homoclinics of the slow--fast symmetric model (see Figure \ref{f-2Dasym-LAH}). As such, the appropriate homoclinic orbit along which to measure the distance function is the concatenation, $\Gamma_0^- \cup \Gamma_0^+$, of the left and right homoclinics to $P_0$. 
Analytically, it is known from multi-pulse Melnikov theory that the Melnikov functions for double-loop homoclinic orbits are given to leading order by the sum of the individual single-loop Melnikov functions, with the contribution from the passage near the saddle being of higher order, see for example \cite{Camassa1998}. 
Thus, the bifurcation equation $D(\lambda) = 0$ for the persistence of the large-amplitude homoclinics to $P_0$ gives
\begin{equation} 
\lambda = \lambda_0^\cup(\delta) = \frac{I_2^+ + I_2^-}{I_0^+ + I_0^-} + \mathcal{O}(\eta).
\Label{eq:lambdaLAHP0}
\end{equation}
The function $\lambda_0^{\cup}(\delta)$ is shown in Figure \ref{f-2Dasym-lambda-p0} (green curve).

The definition of $\mu$ and $\lambda$ (see \eqref{FSasym-rescaling}) generates a linear relation between $p$ and $r$
\[ (\lambda - \mu)(r-1) = \lambda(p-1). \]
This formula yields the slope of the tangent lines to the curves of homoclinic bifurcations at the organizing centers.
Here, $\mu = 1$, and $\lambda$ is positive along both branches of homoclinics. This analysis agrees with the numerical continuation results seen previously. 

\begin{remark}
The unfolding and Melnikov analysis for the homoclinic orbits that emanate from $Q_1$ is similar. 
The main difference is that the hyperbolic saddle of the unperturbed Hamiltonian is located at the equilibrium $P_1$ instead of $P_0$. 
As such, the first step in the analysis is a translation to place $P_1$ at the origin. 
The details are presented in Appendix \ref{sss-Q2homoclinics}.
\end{remark}

\begin{remark}
System \eqref{FSasym-xy'} with $h_\eps$ given to leading order by $h_0=-x$ is a special case of general four-parameter planar vector fields studied in \cite{DG1987,KKR1998}. These four-parameter vector fields arise as part of the unfolding of vector fields, 
$
\dot{x} = y, \,\,
\dot{y} = - x^3 - x^2 y,
$
which have a $\mathbb{Z}_2$-symmetric BT singularity, see \cite{Carr1981,K1985}. In \cite{DG1987}, the four-parameter vector fields are of the form
$
\dot{x} = y,  \,\,
\dot{y} = - (x^3 + ax^2 + bx + c) + (d- x^2)y,
$
where $a, b, c,$ and $d$ are real numbers, and in \cite{KKR1998} the vector field is the same except for one difference,
namely the quadratic term is $xy$. We refer to \cite{DG1987} for an extensive analysis, based on elliptic integrals, of the bifurcation curves in this system, including the Hopf bifurcations, homoclinic bifurcations, saddle-node bifurcations of limit cycles, the effects of the $\mathbb{Z}_2$ symmetry-breaking, and many other bifurcations. 
Also, we observe that with higher--order terms in $h_{\eps}$, the system \eqref{FSasym-xy'} involves higher--order polynomials than those studied in \cite{DG1987,KKR1998}.
We refer to \cite{DRS1991} for the analysis of a related co-dimension three singularity.
\end{remark}

\section{The Maasch--Saltzman Model\label{s-MS}}
In this section, we bring the results from the previous sections together to study the full asymmetric MS model~(\ref{MS-xyz}) for all $q>1$ and $s>0$,
\begin{equation}
\begin{split}
\dot{x} &= - x - y , \\
\dot{y} &= ry - pz + sz^2 - y z^2 , \\
\dot{z} &= - q x - q z .
\end{split}
\Label{MS-xyz'}
\end{equation}
The system~\eqref{MS-xyz'} has two organizing centers
\begin{equation}
Q_0 = \left( \frac{q}{1+q}, \frac{q}{1+q} \right) , \quad
Q_1 = \left( \frac{q}{1+q} + \frac{1}{2} s^2, \frac{q}{1+q} + \frac{1}{4} s^2 \right) .
\Label{MS-Q1Q2}
\end{equation}
We show that, for values of $(p,r)$ near $Q_0$ and all $q>1$, the system \eqref{MS-xyz'} has a family of two-dimensional center manifolds to $P_0$, and there is a critical value $q_c(p,r,s)$ such that the manifolds are at least $C^1$-smooth for all $q > q_c(p,r,s)$.
There is a similar result for values of $(p,r)$ near $Q_1$, and these center manifolds near $Q_1$ coincide with those near $Q_0$.

\subsection{Equilibria and Bifurcations\label{ss-MSEquilStab}}
The equilibria of (\ref{MS-xyz'}) are the same as those in the previous section. 
The trivial state $P_0 = (0,0,0)$ is again a solution for all parameter values. 
If $s^2 + 4(r-p) > 0$, then there are two additional equilibria, $P_1 =  (x_1^*, - x_1^*, - x_1^*)$ and $P_2 =  (x_2^*, - x_2^*, - x_2^*)$, where
\begin{equation}
x_1^* = \tfrac12 [- s + \sqrt{s^2 + 4(r - p)}] , \quad
x_2^* = \tfrac12 [- s - \sqrt{s^2 + 4(r - p)}] ,
\Label{MS-x1*x2*}
\end{equation}
recall (\ref{FSasym-x1*x2*}). Also recall that $x_2^* < x_1^* < 0$ if $p - \tfrac14 s^2 < r < p$, and $x_2^* < 0 < x_1^*$ if $r > p$.

\begin{figure}[h!]
\centering
\includegraphics[width=4in]{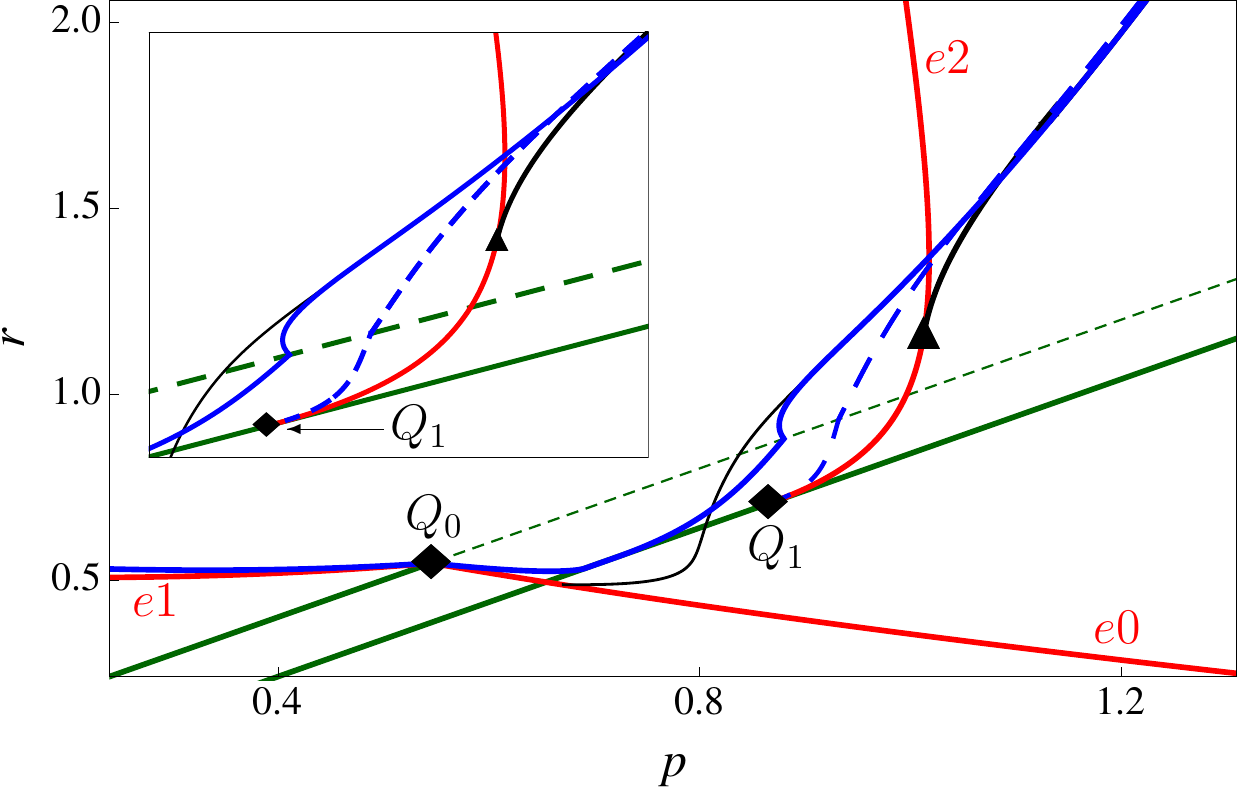}
\caption{Bifurcation structure of \eqref{MS-xyz'} for $q = 1.2$ and $s = 0.8$.
\label{f-MS-stability0}}
\end{figure}

Let $P = (x^*, -x^*, -x^*)$ be any of the equilibria, with $x^* = 0$, $x_1^*$, or $x_2^*$. The characteristic equation of the Jacobian at $P$ is $\lambda^3 + b \lambda^2 + c \lambda  + d = 0$, with
\[ b = 1 + q - r + (x^*)^2 , \quad c = q - (1 + q) r + (1 + q) (x^*)^2 , \quad d = q \left[ p - r + 2 s x^* + 3 (x^*)^2 \right]. \]
By the Routh--Hurwitz conditions, $P$ is (linearly) stable if $b>0$, $c>0$, $d>0$, and $e = bc-d>0$. We analyze these conditions for fixed $q$ and $s$, considering $b$, $c$, $d$, and $e$ as functions of $p$ and $r$. Figure~\ref{f-MS-stability0} illustrates the results for $q = 1.2$ and $s = 0.8$ (the values chosen in~\cite{MaaschSaltzman1990}), using the same color scheme as in Figure \ref{f-2Dasym-BifurcationCurves}.

Let $\Omega$ be the domain in the positive quadrant of the $(p,r)$ plane where all four inequalities are satisfied. The Routh--Hurwitz conditions cease to be satisfied when at least one of the inequalities becomes an equality. Hence, $\partial \Omega$ consists of segments where at least one of the coefficients $b$, $c$, $d$, or $e$ vanishes. 

Consider the points of $\partial \Omega$ where at least one of the coefficients $b$, $c$, $d$, or $e$ vanishes. 
We claim that at any such point,
with possibly finitely many exceptions,
$b>0$, $c>0$, and either $d=0$ and $e>0$
or $d>0$ and $e=0$.
To prove the claim, we assume that at such a point
exactly one of the coefficients $b$, $c$, $d$, or $e$ vanishes,
while the remaining three coefficients are all strictly positive.
This assumption is true generically---that is, 
with possibly finitely many exceptions.
If $b = 0$ or $c=0$, then $e = -d$, 
so $d$ and $e$ cannot both be positive.
Therefore, it must be the case that $b>0$, $c>0$,
and either $d=0$ or $e=0$, as claimed.

If $d=0$, the characteristic polynomial reduces to
$\lambda (\lambda^2 + b \lambda+ c) = 0$,
which yields a simple root at the origin 
and two roots with negative real parts.
If $e=0$, it reduces to
$(\lambda + b)(\lambda^2 + c) = 0$,
which yields a negative real root $-b$
and a conjugate pair of purely imaginary roots $\pm i \sqrt{c}$,
indicating that the equilibrium loses stability
due to a Hopf bifurcation.
The quantity $\sqrt{c}$ is the natural frequency.

\subsubsection{Stability of $P_0$\label{sss-MSEquilStabP0}}
If $x^* = 0$, the coefficients of the characteristic polynomial are $b_0 = 1 + q - r$, $c_0 = q - (1 + q) r$, and $d_0 = q (p - r)$.
The condition $d_0 > 0$ is satisfied if $p > r$.
The condition $e_0 = b_0 c_0-d_0 > 0$ is satisfied if
\begin{equation}
p - r < \frac{1 + q}{q} \left( \frac{q}{1 + q} - r \right)(1 + q - r)  .
\Label{MSsym-P0-RHe}
\end{equation}
Since $p > r$ is a necessary condition for stability, the expression in the right member must be positive---that is, either $0 < r < \frac{q}{1 + q}$ which corresponds to the case where both factors in the right member are positive, or $r > 1 + q$ which corresponds to the case where both factors in the right member are negative. The latter inequality violates the Routh--Hurwitz condition $b_0 > 0$; hence, stability can occur only if $0 < r < \frac{q}{1 + q}$.

The zero-level set of $e_0$ constitutes a curve, labeled e0 in Figure \ref{f-MS-stability0}, of supercritical Hopf bifurcations of $P_0$. The equilibrium $P_0$ loses stability, and a stable limit cycle is created as $(p, r)$ crosses e0 going upward. The curve e0 is a parabola in the $(p, r)$ plane, which is open to the right. The lower branch of the parabola is given by 
\begin{equation} 
r = f_0 (p) = \frac{1}{2} (1+q) - \sqrt{ \frac{1}{4} (q - 1)^2 + \frac{q}{1+q} p} ,
\Label{WhyDoesThisDeserveALabel?}
\end{equation}
where $f_0(p)$ is monotonically decreasing. (We note that the condition $0 < r < \frac{q}{1 + q}$ is satisfied if and only if $\frac{q}{1 + q} < p < 1 + q$.)
Therefore, the Routh-Hurwitz criteria are satisfied in the region enclosed by the $p$ axis, the diagonal, and the curve e0.

In addition, the Jacobian of \eqref{MS-xyz'} at $P_0$ has a zero eigenvalue of geometric multiplicity two at $Q_0$, and a third eigenvalue $\lambda_3=-(1+q - \frac{q}{1+q})$. Hence, $Q_0$ is an organizing center. 

\subsubsection{Stability of $P_1$\label{sss-MSEquilStabP1}}
Recall that $P_1$ exists if and only if $r \ge p - \tfrac14 s^2$. With $x^* = x_1^*$, we have
\begin{equation}
\begin{split}
b_1 &= 1 + q - p + \tfrac12 s^2 - \tfrac12 s \sqrt{s^2 + 4 ( r - p)}] , \\
c_1 &= (1 + q) \left[ \frac{q}{1 + q}  - p + \tfrac12 s^2 - \tfrac12 s \sqrt{s^2 + 4 ( r - p)} \right] , \\
d_1 &= \tfrac12 q \left[s^2 + 4(r - p) - s \sqrt{s^2 + 4(r-p)} \right] .
\end{split}
\Label{MS-P1-RH1}
\end{equation}
The zero-level set of $d_1$ is the diagonal $r=p$. On the diagonal, $P_0$ and $P_1$ exchange stability in a transcritical bifurcation.
The zero-level set, e1, of $e_1 = b_1c_1-d_1$ emerges from $Q_0$ (see Figure \ref{f-MS-stability0}) and corresponds to a curve of supercritical Hopf bifurcations of $P_1$. Therefore, the Routh-Hurwitz criteria are satisfied in the region enclosed by the $r$ axis, the diagonal, and the curve e1. 

In addition, the Jacobian of \eqref{MS-xyz'} at $P_1$ has a zero eigenvalue of geometric multiplicity two at $Q_1$, and a third eigenvalue $\lambda_3=-(1+q - \frac{q}{1+q})$. Hence, $Q_1$ is an organizing center.

\subsubsection{Stability of $P_2$\label{sss-MSEquilStabP2}}
Recall that $P_2$, like $P_1$, exists if and only if $r \ge p - \tfrac14 s^2$. With $x^* = x_2^*$, we have
\begin{equation}
\begin{split}
b_2 &= 1 + q - p + \tfrac12 s^2 + \tfrac12 s \sqrt{s^2 + 4 ( r - p)}]   , \\
c_2 &= (1 + q) \left[ \frac{q}{1 + q}  - p + \tfrac12 s^2 + \tfrac12 s \sqrt{s^2 + 4 ( r - p)} \right] , \\
d_2 &= \tfrac12 q \left[ s^2 + 4(r - p) + s \sqrt{s^2 + 4(r-p)} \right] .
\end{split}
\Label{MS-P2-RH}
\end{equation}
The zero-level set of $d_2$ is the shifted diagonal $r = p - \tfrac14 s^2$. On the shifted diagonal, $P_1$ and $P_2$ are created in a saddle-node bifurcation. 
The zero-level set, e2, of $e_2 = b_2c_2-d_2$ is the continuation of e1 (after a gap between $Q_0$ and $Q_1$), and corresponds to a curve of subcritical Hopf bifurcations of $P_2$.
At $Q_1$, the zero-level set of $e_2$ is tangent to the shifted diagonal. 
Therefore, the Routh-Hurwitz criteria are satisfied in the region enclosed by the axes, the shifted diagonal, and e2. 

Along the Hopf bifurcation curve e2, there is a Bautin bifurcation point, marked by the black triangle in Figure \ref{f-MS-stability0}. The Hopf bifurcations are subcritical below the Bautin point and supercritical above it. Also, a (black) branch of saddle-node bifurcations of limit cycles emerges  from the Bautin point. 
For each fixed $s\ge 0$, the Bautin point is close to $Q_1$ for $q$ close to one, and then it moves up along the curve e2 as $q$ increases.

\begin{remark}
In the symmetric limit ($s \to 0$), the curves e1 and e2 collapse to the same parabola in the $(p,r)$ plane, given by 
\begin{equation}
p = f_{1,2}(r) = \frac{1}{2}\left( 1+q-\frac{q}{1+q} \right)-\sqrt{\frac{1}{4} \left( 1+q-\frac{q}{1+q} \right)^2 -q+\frac{2q}{1+q}r}.
\Label{e1e2}
\end{equation}
That is, the curves e1 and e2 are different $s$-unfoldings of the parabola \eqref{e1e2}.  
\end{remark}

\subsubsection{Global Bifurcations\label{sss-MSglobalbifurcations}}
The curves of global bifurcations that emanate from the organizing centers $Q_0$ and $Q_1$ are shown in Figure \ref{f-MS-stability0}. The blue curve is the curve of homoclinic bifurcations of \eqref{MS-xyz'}. The system \eqref{MS-xyz'} possesses the same six types of single-loop and double-loop (large-amplitude) homoclinic orbits to the saddles $P_0$ and $P_1$ as the system \eqref{FSasym-xy} on the slow manifold $\mathcal{M}_\eps$. They lie in regions I and III (recall Figure \ref{f-2Dasym-regions}),
where the lower boundaries of these regions are now given by the Hopf bifurcation curves e1 and e0, respectively. 

The black curves in Figure \ref{f-MS-stability0} indicate curves of saddle-node bifurcations of limit cycles. Along these black curves, a pair of limit cycles of \eqref{MS-xyz'} of opposite stability merge and annihilate each other. The existence of these global bifurcation curves follows directly from the existence of center manifolds (see the next subsection) and the presence of BT points of the reduced equations on the center manifolds.

\begin{remark}
Figure~\ref{f-MS-SymmetryBreaking} shows how the bifurcation structure of \eqref{MS-xyz'} collapses as $s \to 0$. 
The shifted diagonal of saddle-node bifurcations of $P_1$ and $P_2$ converges to the diagonal of transcritical bifurcations, resulting in a diagonal line of pitchfork bifurcations in the symmetric limit.
The organizing center $Q_1$ collapses onto the organizing center $Q_0$, leaving the $\mathbb{Z}_2$-symmetric BT point, $Q$, seen in the symmetric MS model~\eqref{MSsym-xyz-1}. Consequently, the curves e1 and e2 of Hopf bifurcations of $P_1$ and $P_2$ merge to the curve \eqref{e1e2}. Similarly, the curves of homoclinic bifurcations merge to a single curve of homoclinics, and the curves of saddle-node bifurcations of limit cycles coalesce and become a single curve. 

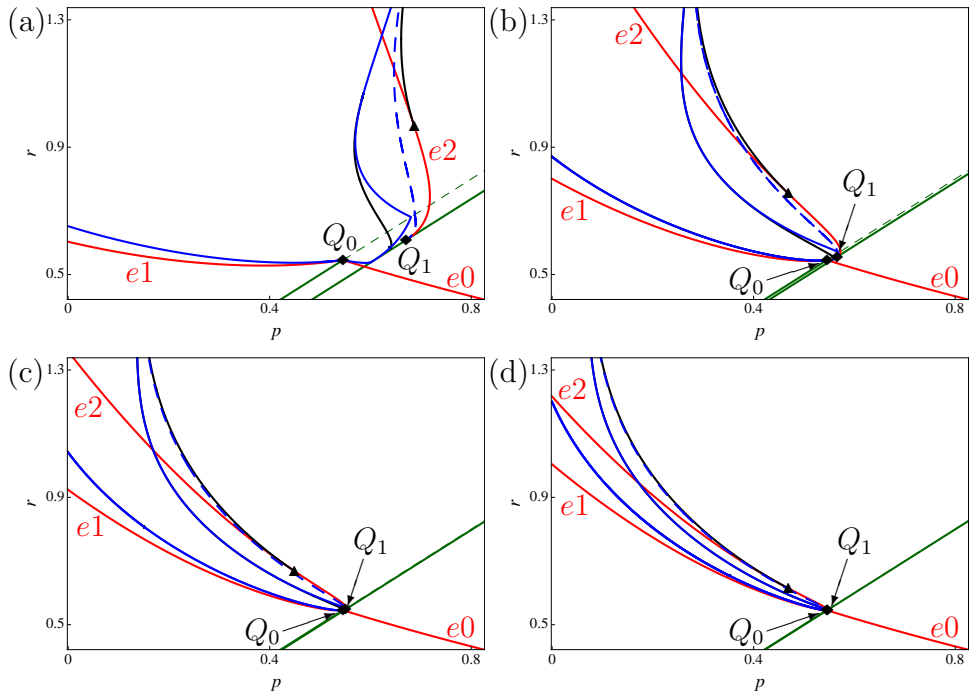
\begin{figure}[h!]
\centering
\includegraphics[width=5in]{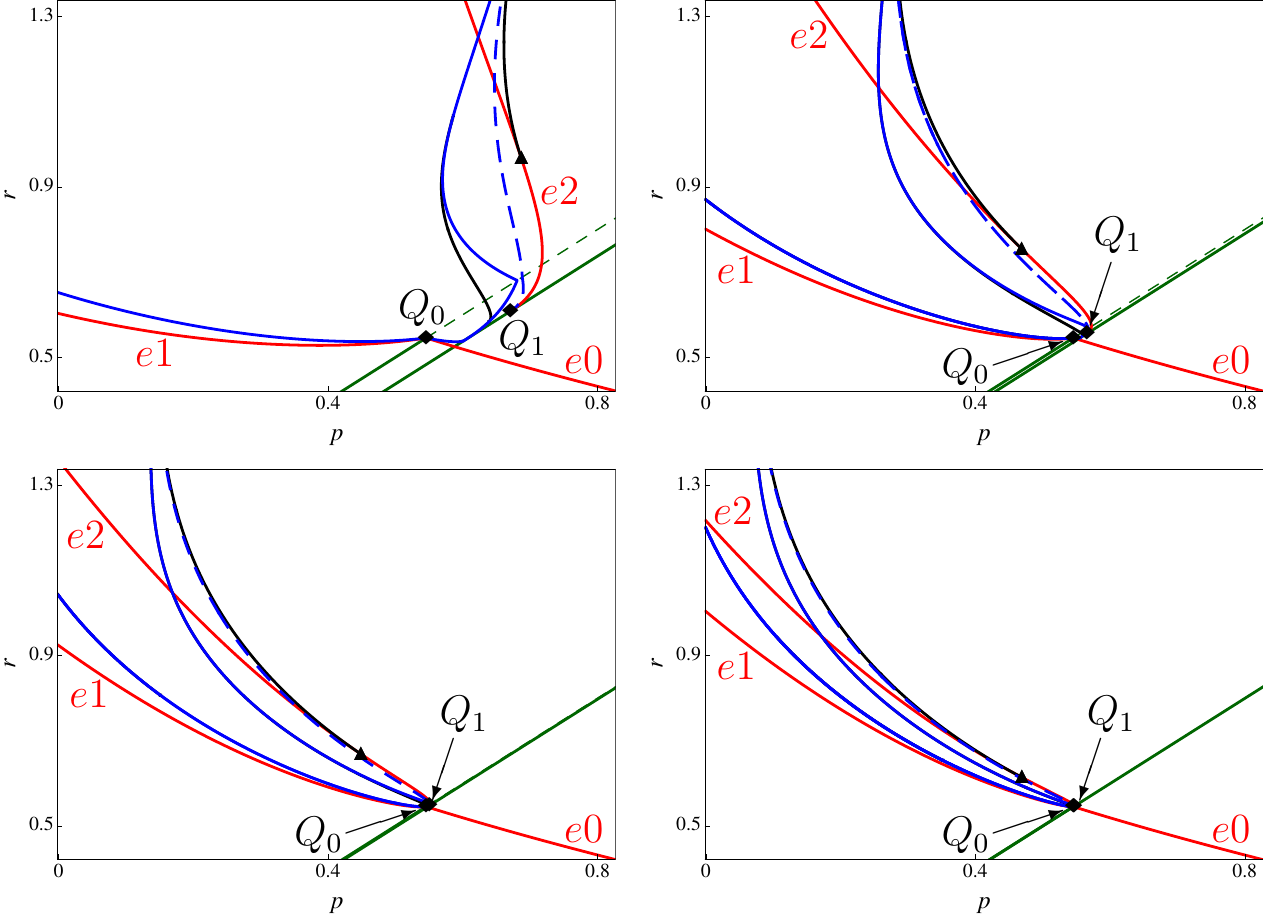}
\put(-368,250){(a)}
\put(-186,250){(b)}
\put(-368,117){(c)}
\put(-186,117){(d)}
\caption{Bifurcation curves of \eqref{MS-xyz'} with $q = 1.2$, (a) $s = 0.5$, (b) $s=0.2$, (c) $s=0.1$ and (d) $s=0.05$.
\label{f-MS-SymmetryBreaking}}
\end{figure}

\end{remark}



\subsection{Center Manifolds\label{a4-ss-CenterManifolds}}
In this section, we establish the existence of the center manifolds associated to the equilibria $P_0$ and $P_1$ in the system~(\ref{MS-xyz'}) for parameter values near the organizing centers $Q_0$ and~$Q_1$. We follow the general approach described in~\cite[\S 3.2]{Guckenheimer1983}. The derivation of the center manifolds of $P_0$ for parameters near $Q_0$ is given in detail. Since $P_1$ lies on the center manifolds to $P_0$, the center manifolds attached to $P_1$ are in fact the same as those attached to $P_0$ to all orders.

As a first step to establish the center manifolds of $P_0$, we transfer $Q_0$ to the origin in parameter space by introducing $(\tilde{p},\tilde{r}) = (p - \frac{q}{1 + q}, r - \frac{q}{1 + q})$. Thus, \eqref{MS-xyz'} transforms to 
\begin{equation}
\begin{pmatrix} \dot{x} \\ \dot{y} \\ \dot{z} \end{pmatrix}
= A \begin{pmatrix} x \\ y \\ z \end{pmatrix}
+ \begin{pmatrix}
0 \\ n(x,y,z,\tilde{p},\tilde{r},s) \\ 0
\end{pmatrix} ,
\Label{MS-vector-xyz}
\end{equation}
where
\begin{equation}
A = \begin{pmatrix}  -1 & -1 & 0 \\   0 & \frac{q}{1 + q} & - \frac{q}{1 + q} \\   -q & 0 & -q  \end{pmatrix} 
\quad \text{ and } \quad
n(x,y,z,\tilde{p},\tilde{r},s) = \tilde{r}y - \tilde{p}z + s z^2 - y z^2. 
\Label{MSsym-A,n}
\end{equation}
The matrix $A$ has eigenvalues $\lambda_1 = 0$, $\lambda_2 = 0$, and $\lambda_3 = - \left( 1+q-\frac{q}{1+q} \right)$.
%
%
Next, we change coordinates to reduce $A$ to its Jordan normal form,
\begin{equation}
J = F^{-1} A F = \begin{pmatrix} 0 & 1 & 0 \\ 0 & 0 & 0 \\ 0 & 0 & \lambda_3 \end{pmatrix} .
\end{equation}
The columns of $F$ are the (generalized) eigenvectors of~$A$,
\begin{equation}
F = \begin{pmatrix} 1 & 1 & \frac{1}{q (1 + q)} \\ -1 & -2 & \frac{q}{(1 + q)^2} \\ -1 & \frac{1 - q}{q} & 1 \end{pmatrix}.
\end{equation}

If $(x, y, z)$ satisfies~(\ref{MS-vector-xyz}), then $(u, v, w)^T = F^{-1} (x, y, z)^T$ satisfies 
\begin{equation}
\begin{pmatrix} \dot{u} \\ \dot{v} \\ \dot{w} \end{pmatrix}
= J \begin{pmatrix} u \\ v \\ w \end{pmatrix} + \frac{1 + q}{\ell^2} \, n \begin{pmatrix} q \ell - 1 \\ - q \ell \\ q (1 + q) \end{pmatrix} ,
\Label{MS-uvw}
\end{equation}
where we note that the variables $u$ and $v$ here are distinct from those used in the unfolding analysis of the previous section, we have introduced the abbreviation $\ell = 1 + q + q^2$, and $n(u,v,w,\tilde{p},\tilde{r},s)$ is 
\begin{equation}	\label{eq:nonlinearuvw}
\begin{split}
n(u,v,w,\tilde{p},\tilde{r},s) &=
\tilde{r}
\left(  - u - 2v + \frac{q}{(1+q)^2} w \right)
- \tilde{p}
\left(  - u + \frac{1-q}{q} v + w \right) \\
& \hspace{2em} - \left( - s - u - 2v + \frac{q}{(1+q)^2} w \right)
\left(  - u + \frac{1-q}{q} v + w \right)^2 .
\end{split}
\end{equation}

By standard center manifold theory \cite{Carr1981,Guckenheimer1983,Kuznetsov2004}, there exists a family of two-dimensional center manifolds, $W^c(0)$, which for any $k>0$ are given by
\begin{equation}
W^c(0) = \{ (u, v, w) : w = h(u, v, \tilde{p}, \tilde{r}), \, h \in C^k \}.
\Label{MSsym-CenterManifold}
\end{equation}
The center manifolds are not unique, but they are $C^k$ equivalent. The function $h$ satisfies the invariance equation,
\begin{equation}
\frac{\partial h}{\partial u} \dot{u} + \frac{\partial h}{\partial v} \dot{v} = \lambda_3 h + \frac{q(1 + q)^2}{\ell^2} \, n,
\Label{MS-invariance}
\end{equation}
with $h=0$, $\frac{\partial h}{\partial u}=0$, and $\frac{\partial h}{\partial v}=0$ at $(0,0,\tilde{p},\tilde{r})$. Also, $h$ may be represented by a series of the form
\begin{equation}
h = h_0 + h_1 + h_2 + h_3 + \cdots ,
\Label{MS-h}
\end{equation}
where $h_0$ is a constant and $h_i$ ($i = 1, 2, \ldots$) is a homogeneous polynomial function of degree $i$ of the variables $u$, $v$, $\tilde{p},\tilde{r}$, and the coefficients in these polynomials depend on $q$ and $s$.
The first two terms vanish identically since the center manifolds are tangent to the center subspace at the origin. The expressions for $h_2$ and $h_3$ are given in Appendix~\ref{app:CM}.

\begin{remark}
It is useful to compare the expression~(\ref{MS-h}) for $h(u,v,\tilde{p},\tilde{r})$ with the expression~(\ref{hepsasymm}) for $h_{\eps}(x,y,p,r)$ for the slow--fast system. 
Using the transformation $(x, y, z)^T = F(u,v,w)^T$ and the slow manifold expansion~(\ref{hepsasymm}) with terms up to and including $\mathcal{O}(\eps^3)$, one finds that 
\begin{align*} 
\left| u+(1-\eps)v + h_{\eps} \left( x (u, v, h (u, v, \tilde{p}, \tilde{r})), y (u, v, h(u, v, \tilde{p}, \tilde{r}) ), \tilde{p}, \tilde{r} \right) - h(u, v, \tilde{p}, \tilde{r}) \right| \\ 
= \mathcal{O} \left( \eps^4, \eps^3 \left( u+v+p+r \right)^4 \right),
\end{align*}
where the quantity $u+(1-\eps)v+h_{\eps}$ is the representation of the slow manifold in the $(u,v,w)$ coordinates. 
Thus, the slow manifold and the center manifold are $\mathcal{O}(\eps^4)$ close as $\eps \to 0$ ($q \to \infty$), which is an important consistency check.
This shows that the center manifolds naturally continue to the slow manifolds from the regime~$q \gg 1$. 
\end{remark}

\begin{figure}[h!]
\centering
\includegraphics[width=3in]{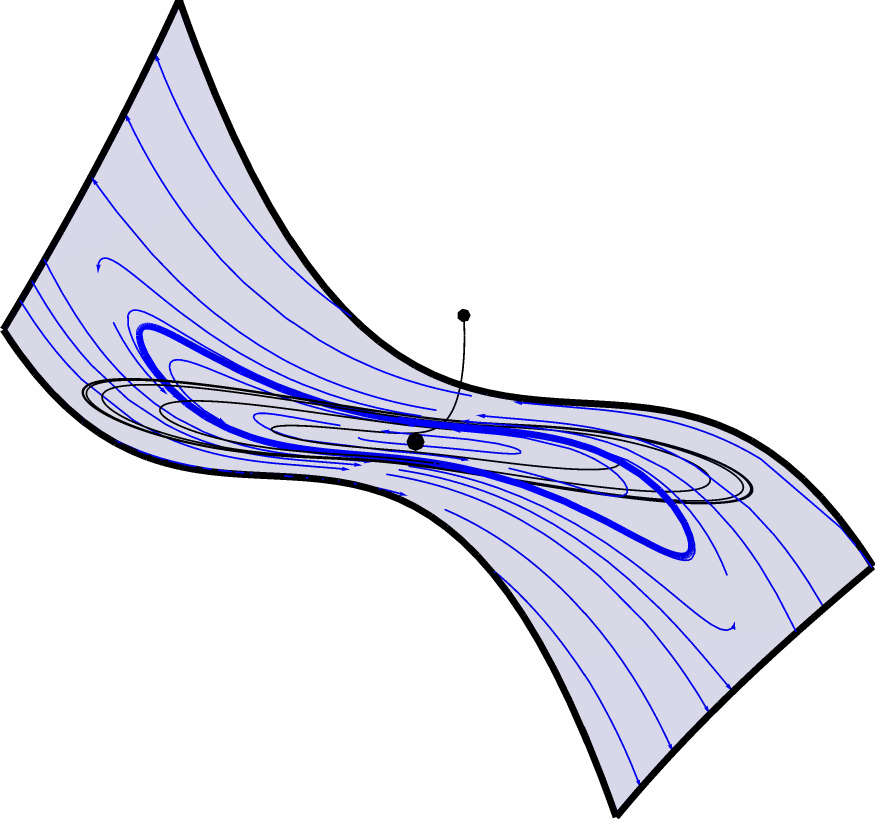}
\caption{The local center manifold (blue surface) for $s=0$,  $q=1.2$ and $(\tilde{p}, \tilde{r}) = (0.15, 0.1)$, the flow on the center manifold (blue streamlines), the stable limit cycle (thick blue) on the center manifold, and a trajectory (black curve) of~\eqref{MS-vector-xyz} with initial condition (black marker) away from $W^c (0)$. Note that the stable limit cycle  deviates from the actual trajectory due to the fact that the center manifold is approximated up to and including cubic terms.
\label{f-MSsym-redeqonWc}} 
\end{figure}

On the center manifold, the dynamics are governed by
\begin{equation}
\begin{split}
\dot u &= v + \frac{(1+q)(q\ell -1)}{\ell^2} \, n (u, v, h (u, v, \tilde{p}, \tilde{r}), \tilde{p}, \tilde{r},s) , \\
\dot v &= \frac{-q (1+q)}{\ell} \, n (u, v, h (u, v, \tilde{p}, \tilde{r}), \tilde{p}, \tilde{r},s) ,
\end{split}
\Label{MS-redeqonWc}
\end{equation}
where $n$ is given by~\eqref{eq:nonlinearuvw}. The center manifold $W^c(0)$ and the reduced equations \eqref{MS-redeqonWc} on it are illustrated in Figure~\ref{f-MSsym-redeqonWc}, with $s=0$ and $(p,r)$ chosen such that $P_0$ is an unstable focus enclosed by a stable limit cycle. Figure~\ref{f-MSsym-redeqonWc} also shows how a solution through a nearby initial condition rapidly approaches $W^c(0)$ and then winds towards the stable limit cycle.


The system \eqref{MS-redeqonWc} has a BT at $Q_0$ corresponding to the equilibrium $P_0$, and the unfolding of the BT point is similar to that in Section~\ref{ss-BTanalysis} (details not shown). Also, as shown in Figures \ref{f-MS-stability0} and \ref{f-MS-SymmetryBreaking}, the bifurcation curves that emanate from the BT point on $W^c(0)$ are of the same type as those that emanate from $Q_0$ on the slow manifold. 

The smoothness of the center manifolds is determined by the Lyapunov type numbers that measure the growth rates of solutions in the directions tangential to the center manifolds and the growth rates normal to the manifolds \cite{Fenichel1971}. For each $(p,r,s)$, there is a critical value of $q$, which we label $q_c(p,r,s)$, given by the largest value of $q > 1$ for which both Lyapunov type numbers are less than one. The center manifolds are at least $C^1$-smooth for all $q>q_c(p,r,s)$.

Numerically, the Lyapunov type numbers can be approximated from the eigenvalues at the equilibria of \eqref{MS-xyz'} and the Floquet multipliers of the limit cycles of \eqref{MS-xyz'}. Let $\sigma_{ij} (p,r,s)$, $j=1,2$, denote the eigenvalues at the equilibrium $P_i$, $i=0,1,2$, corresponding to the two eigenvectors in the tangent plane to the manifolds. For $s=0$ and for each $(p,r)$ in an array of values above the diagonal, we calculated $M(p,r) = \max_{i=0,1,2; j=1,2} \vert \operatorname{Re}(\sigma_{i,j} (p,r,0)) \vert$. This quantity is a function of $q$, and we determined the smallest value of $q$ for which $M(p,r) < |\lambda_3|$, where we recall that $\lambda_3=-(1+q-(q/1+q))$ is the transverse eigenvalue at $Q_0$. This smallest value marks the value at which the tangential growth rate at the equilibria on $W^c(0)$ first equals the normal growth rate. This enables a simple numerical approximation of $q_c$, since for all $q$ greater than this value, the tangential growth rate is less than the normal growth rate, as approximated by $\lambda_3$, and hence shows that the reduced equations \eqref{MS-redeqonWc} are an approximation to \eqref{MS-xyz'}. 
Numerically, for $(p,r) \in \left\{ 0 < p \leq 2, 0 < r \leq \tfrac{3}{2} \right\}$, we find that $q_c(p,r,0)$ is less than or equal to one.  For $p$ fixed at $1$, we find that $q_c(p,r,0)$ increases monotonically with $r$ (with $q_c(1,2,0) \approx 1.4$). Similarly, $q_c$ increases monotonically with $r$ for other fixed values of $p$. 
We also computed the Floquet multipliers of the limit cycles on $W^c(0)$ for sample points $(p,r)$ in the same array, and verified that the maximal tangential growth rates associated to the limit cycles are less than $\vert \lambda_3 \vert$ for $q>q_c(p,r)$. 

\begin{remark}
The cubic approximation of $W^c(0)$ is valid in a neighborhood of $P_0$. Any extraneous equilibria generated by it lie outside the domain of validity.
\end{remark}

\section{Discussion}     \label{s-Discussion}

In this article, we presented a dynamical systems analysis of the rich behavior exhibited by~(\ref{MS-xyz}), which was proposed by Maasch and Saltzman in~\cite{MaaschSaltzman1990} in their study of the glacial cycles observed in the climate record of the Pleistocene Epoch. We identified the regimes in which the MS model~(\ref{MS-xyz}) exhibits limit cycles, and we determined the locations of the various bifurcation curves along which the limit cycles are created and disappear. 

Central to understanding the limit cycles and their bifurcations is the result that the long-term system dynamics of the third-order MS model~(\ref{MS-xyz})
occur on (and near) two-dimensional invariant manifolds
for most values of the parameter~$q$,
which measures the characteristic time scales
for the total global ice mass
and the volume of North Atlantic Deep Water.
First, by considering the regime
in which $q$ is asymptotically large~($q \gg 1$),
we showed in Sections~2 and~3 that
the model possesses two-dimensional invariant slow manifolds.
All initial conditions relax quickly to these manifolds,
and the solutions approach the stable equilibria and stable limit cycles
along the manifolds.
These slow manifold results
generalize the earlier analysis
presented in~\cite{EKKV2017}
for the second-order system
obtained by formally setting $q=\infty$
and $s=0$ in~(\ref{MS-xyz}).

Second, by considering finite---but not large---values of $q$,
we showed in Section~4 that
the model has a family
of two-dimensional invariant center manifolds
for all finite values of $q>1$.
These center manifolds
also contain the equilibria and limit cycles,
attract all nearby initial conditions,
and govern how solutions approach the stable states.
Moreover, the center manifolds are at least $C^1$-smooth for all $q$ greater than a critical value $q_c$. 
We note that the center manifolds
for finite $q$
smoothly transition
to the slow manifolds
as $q$ becomes asymptotically large,
as is consistent
with the theory that slow manifolds
may be viewed as special types
of center manifolds~\cite{Carr1981, F1978}.

On both the slow and center invariant manifolds,
there are Bogdanov-Takens points
which act as organizing centers of the dynamics.
These were first studied
in the symmetric case
in which the parameter~$s$ was set to zero.
There is a single $\mathbb{Z}_2$-symmetric BT point
from which the three main types of bifurcation curves emanate:
the Hopf bifurcation curves
along which the stable and unstable limit cycles
are created from the equilibria,
the homoclinic bifurcation curve
along which the pair of small-amplitude unstable limit cycles
coalesce to form one large-amplitude unstable limit cycle,
and the curve of saddle-node bifurcations of limit cycles
along which the large-amplitude stable and unstable limit cycles disappear.
Then, it was shown that
there is a symmetry-breaking
which creates two non-degenerate BT points
in the physically-relevant case of~$s>0$,
where the limit cycles
of~(\ref{MS-xyz})
exhibit slow glaciation and rapid deglaciation.
The bifurcation curves emanating from
these two organizing centers
are similar to,
but more complex than, those in the symmetric case.

\begin{figure}[h!]
\centering
\includegraphics[width=4in]{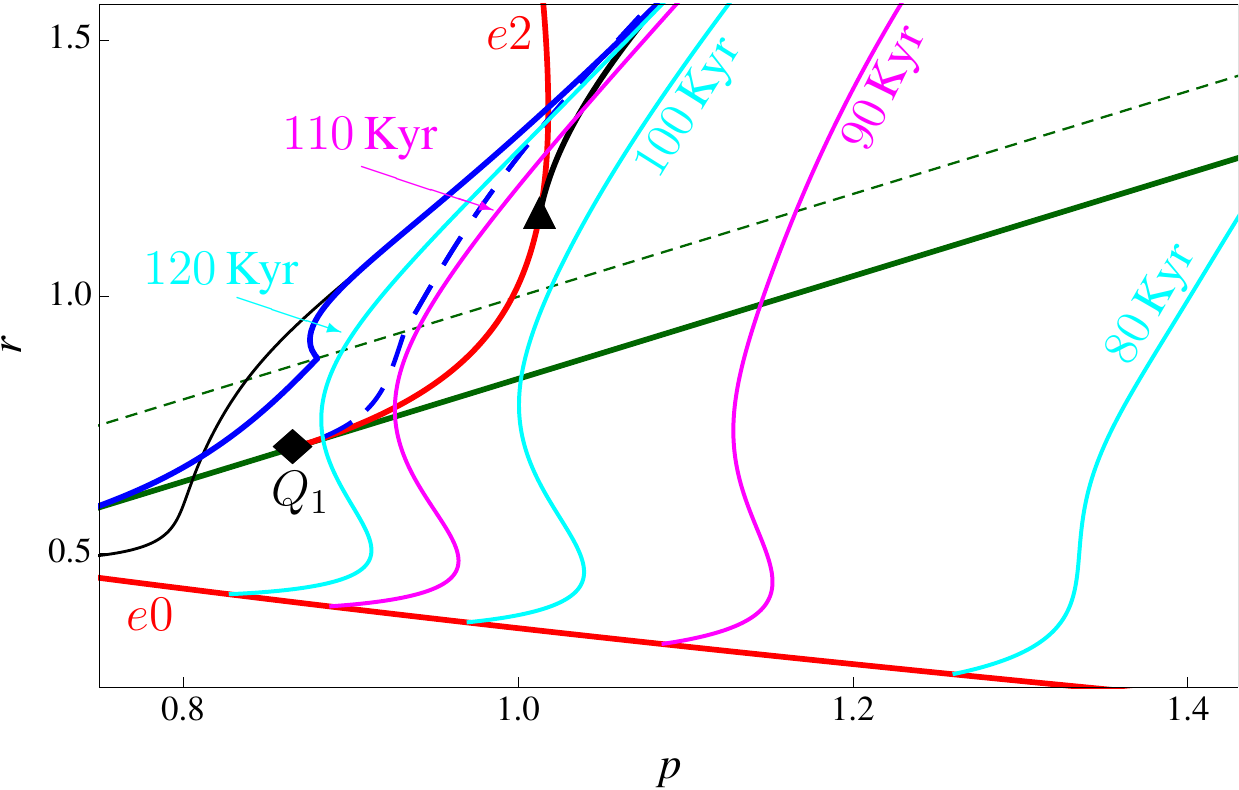}
\caption{Isoperiod curves (alternating pink and cyan) for~\eqref{f-2Dasym-isoperiods}
with $s = 0.8$ and $q=1.2$.
\label{f-2Dasym-isoperiods}}
\end{figure}

A summary of the local and global bifurcation curves
in which the limit cycles are created and disappear
is given by Figures~\ref{f-MS-stability0}
and~\ref{f-MS-SymmetryBreaking}.
Stable limit cycles
of the full MS model~(\ref{MS-xyz})
are found in the region bounded
by the curve e0 of supercritical Hopf bifurcations of $P_0$
and the curve of homoclinic bifurcations
(solid blue branch emanating
from~$Q_0$ to the right).
This same region is also shown
in Figure~\ref{f-2Dasym-isoperiods}
now with the isoperiod curves
(alternating pink and cyan curves)
for the stable limit cycles
super-imposed.
Maasch and Saltzman~\cite{MaaschSaltzman1990}
focused on the parameter region
near the isoperiod curve corresponding to 100\,Kyr cycles.

Besides being of intrinsic interest,
the results presented herein
for the internal dynamics of the MS model~(\ref{MS-xyz})
will be useful for studying
the effects of orbital (Milankovitch) forcing and
slow variation of the system parameters in~(\ref{MS-xyz}).
In~\cite{MaaschSaltzman1990},
Maasch and Saltzman showed
that the model can be tuned
to exhibit the 40\,Kyr cycles of the early Pleistocene
under orbital (Milankovitch) forcing, and
that slow passage of the system parameters
$p$ and $r$ through Hopf bifurcations
offers a mechanism to understand the mid-Pleistocene transition.
We are presently pursuing these findings.

\begin{appendix}
\section{Persistence of the Homoclinics to $P_1$\label{sss-Q2homoclinics}}
In this appendix, we present the Melnikov analysis for the existence and persistence of the single- and double-loop homoclinics to the saddle $P_1$ for the slow--fast asymmetric system. We consider~\eqref{FSasym-uv} in the regions~III and~IV, where $-\tfrac{1}{4}\delta^2 < \mu < 0$. Without loss of generality, we set $\mu=-1$. Then \eqref{FSasym-uv} possesses three equilibria, located at
\[ (u_0,v_0) = (0,0), \quad (u_{1,2}^*,v_{1,2}^*) = \left( \frac{1}{2} \left( -\delta \pm \sqrt{\delta^2-4} \right), 0 \right). \]
In regions III and IV, the $x$-coordinates of $P_1$ and $P_2$ satisfy $x_2^* < x_1^* <0$.  Thus, we identify $(u_1^*,v_1^*)$ with $P_1$ and $(u_2^*,v_2^*)$ with $P_2$. We also note that $\delta > 2$ in this region (for $\mu=-1$). 

First, we translate $P_1$ to the origin via the coordinate transformation 
\[ u = u_1^* +\overline{u}, \quad v = \overline{v}. \]
After dropping the overlines, the translated system is
\begin{equation} 
\begin{split}
\dot{u} &= v, \\
\dot{v} &= \nu u - \kappa u^2 - u^3 + \eta \left( \lambda - (u_1^*)^2 - 2u_1^* u - u^2\right)v, 
\end{split}
\Label{eq:P1uvsystem}
\end{equation}
where 
%
\[ \nu = -u_1^* \sqrt{\delta^2-4} >0 \quad \text{ and } \quad  \kappa = \delta+3u_1^*. \]
The Hamiltonian of the unperturbed version ($\eta = 0$) of~\eqref{eq:P1uvsystem} is
$H(u,v) = \tfrac12 v^2 - \tfrac12 \nu u^2 + \tfrac13 \kappa u^3 + \tfrac14 u^4$,
and the level set $H=0$ corresponds to a pair of homoclinic orbits, $\Gamma_1^\pm$,
\begin{equation} \label{eq:Gamma0_P1}
\Gamma_1^\pm 
= \left\{ t \mapsto \left( u_1^\pm (t), v_1^\pm (t) \right) 
=
\left( 
\frac{\pm 3\alpha \nu}{\cosh (\sqrt{\nu} t) \pm \alpha \kappa}, 
\frac{\mp 3\alpha \nu \sqrt{\nu} \sinh(\sqrt{\nu} t)}{\left( \cosh(\sqrt{\nu}t) \pm \alpha \kappa \right)^2} 
\right)
\right\} ,
\end{equation}
where $\alpha = \left( \kappa^2+\tfrac{9}{2}\nu \right)^{-1/2}$. 

\begin{figure}[h!]
\centering
\includegraphics[width=3.5in]{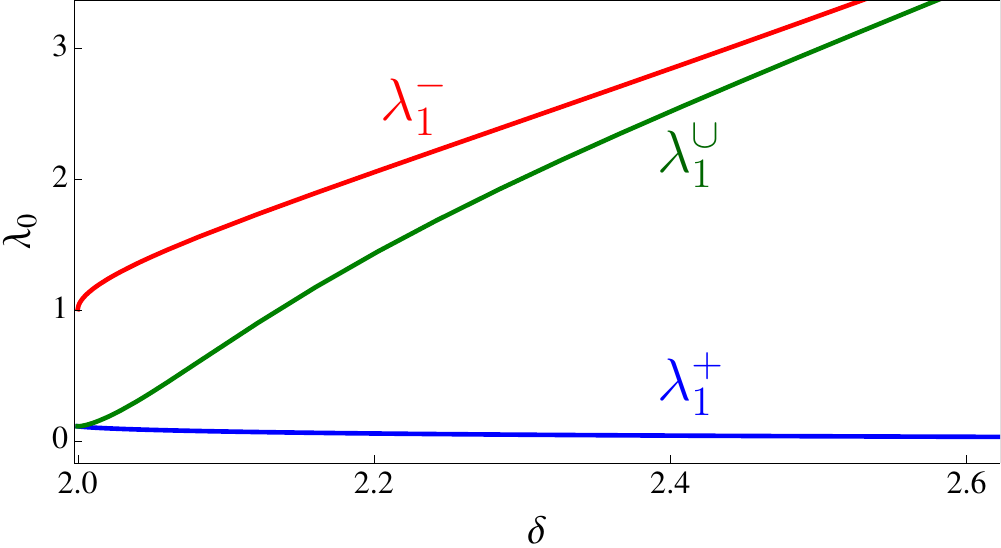}
\caption{Simple zeros $\lambda_1^{\pm}$ and $\lambda_1^\cup$ (\eqref{eq:lambdapm_P1} and \eqref{eq:lambdaLAH_P1}) of the distance function $D(\lambda)$, see \eqref{eq:D_P1}, for $\mu = -1$, corresponding to the right, left, and large-amplitude homoclinics to $P_1$. Recall that $\delta^2 + 4\mu \geq 0$, i.e., $\delta \geq 2$.}
\label{f-2Dasym-lambda-P1}
\end{figure}

The splitting distance is obtained as an asymptotic expansion in $\eta$,
\begin{equation} 
D(\lambda) = \eta \int_{-\infty}^{\infty} \left. \left( \nabla H \cdot \begin{pmatrix} 0 \\ \left( \lambda - (u_1^*)^2 - 2u_1^* u - u^2 \right) v \end{pmatrix} \right) \right|_{\Gamma_1^{\pm}} \, dt + \mathcal{O}(\eta^2). 
\label{eq:D_P1}
\end{equation}
The bifurcation equation $D(\lambda) = 0$ gives the set of $\lambda$ for which the homoclinic orbits persist under perturbations. Thus, we have 
\begin{equation} 
\lambda = \lambda_1^\pm = (u_1^*)^2 + \frac{2u_1^* I_1^\pm + I_2^\pm}{I_0^\pm},
\label{eq:lambdapm_P1}
\end{equation}
where the integrals $I_0^\pm, I_1^\pm,$ and $I_2^\pm$ are defined by
\begin{equation*}
\begin{split}
I_0^\pm &= \int_{-\infty}^{\infty} \left( v_1^\pm \right)^2\, dt = \frac{2}{3}\nu \sqrt{\nu} \left( 3\phi^2-1 \right)+4\nu \sqrt{\nu}\, \psi \phi^2 \tan^{-1} \left( \psi \mp \phi \right), \\
I_1^\pm &= \int_{-\infty}^{\infty} u_1^\pm \, \left( v_1^\pm \right)^2\, dt = -\sqrt{2} \nu^2 \left\{ \psi \left( \frac{5}{2} \phi^2-\frac{1}{3} \right) + \left( 5\psi^4+6 \psi^2+\frac{1}{3} \right) \tan^{-1} (\psi \mp \phi) \right\}, \\
I_2^\pm &= \int_{-\infty}^{\infty} \left( u_1^\pm  \, v_1^\pm \right)^2\, dt = \nu^2 \sqrt{\nu} \left\{ \frac{16}{15} + \frac{23}{3} \psi^2 + 7\psi^4+2\psi \phi^2 \left( 3\phi^2+4\psi^2 \right) \tan^{-1} \left( \psi \mp \phi \right) \right\},
\end{split}
\end{equation*}
where now $\psi = \sqrt{ \frac{2\kappa^2}{9\nu} }$ and $\phi = \sqrt{ \frac{2\kappa^2+9\nu}{9\nu} }$ (and $\phi^2 - \psi^2 = 1$). Figure \ref{f-2Dasym-lambda-P1} shows the simple zeros of $D(\lambda)$.

This demonstrates analytically the existence of the curves of homoclinic bifurcations to~$P_1$. 
Here, $\lambda_1^+$ corresponds to the right-homoclinics (Figure \ref{f-2Dasym-Homoclinics}(b)) which lie along the (flat) branch emanating from $Q_0$ to the right and 
touches the shifted diagonal, and $\lambda_1^-$ corresponds to the left-homoclinics (Figure \ref{f-2Dasym-Homoclinics}(e)) which lie along the dashed branch
that emanates from $Q_1$ to the right, below the diagonal (cf. Figure \ref{f-2Dasym-regions}). 

The large-amplitude homoclinics to $P_1$ for $\mu<0$ are obtained by measuring the splitting distance along the concatenation, $\Gamma_1^- \cup \Gamma_1^+$, of the left and right homoclinics. In this case, the bifurcation equation $D(\lambda) = 0$ yields
\begin{equation} 
\lambda = \lambda_1^\cup (\delta) = (u_1^*)^2+\frac{2u_1^* \left( I_1^+ + I_1^- \right) + (I_2^++I_2^-)}{I_0^++I_0^-} + \mathcal{O}(\eta),
\label{eq:lambdaLAH_P1}
\end{equation}
as the values of $\lambda$ for which the large-amplitude homoclinics to $P_1$ persist. The function $\lambda_1^\cup(\delta)$ is shown in green in Figure \ref{f-2Dasym-lambda-P1}.

The definition~\eqref{FSasym-rescaling} of $\mu$ and $\lambda$ generates a linear relation between $p$ and $r$, namely
\[ (\lambda - \mu) \left( r- (1+ \tfrac14 {s^2}) \right) = \lambda \left( p - (1+\tfrac12 {s^2}) \right) + \tfrac14 {s^2} (\lambda + \mu). \]
This formula yields the slope of the tangent line to the curves of homoclinic bifurcations at~$Q_1$.
With $\mu = -1$, this analysis also agrees with the numerical continuation results. 

\section{Center Manifold Reduction\label{app:CM}}
%
In this appendix, we present the terms in the series~(\ref{MS-h}) for the function~$h$, whose graph represents the center manifold~$W^c (0)$ of the Maasch--Saltzman model~(\ref{MS-xyz'}). Recall that $h = h_2 + h_3 + \cdots$, where~$h_2$ and $h_3$ are homogeneous polynomials of degree 2 and 3, respectively, and the coefficients depend on $q$ and $s$,
\begin{equation}
\begin{split}
h_2 &= b_1 u^2 + b_2 uv + b_3 v^2
+ b_4 \tilde{p}u + b_5 \tilde{p}v
+ b_6 \tilde{r}u + b_7 \tilde{r}v
+ b_8 \tilde{p}^2
+ b_9 \tilde{p} \tilde{r}
+ b_{10} \tilde{r}^2 \\
h_3 &= c_1 u^3 + c_2 u^2 v + c_3 uv^2 + c_4 v^3
+ c_5 \tilde{p}u^2 + c_6 \tilde{p}uv + c_7 \tilde{p}v^2
+ c_8 \tilde{r}u^2
+ c_9 \tilde{r}uv
+ c_{10} \tilde{r}v^2 \\
&\hspace{1em}
+ c_{11} \tilde{p}^2 u + c_{12}\tilde{p}^2 v
+ c_{13} \tilde{p}\tilde{r}u + c_{14} \tilde{p}\tilde{r}v
+ c_{15} \tilde{r}^2 u + c_{16} \tilde{r}^2 v
+ c_{17} \tilde{p}^3
+ c_{18} \tilde{p}^2 \tilde{r} \\
&\hspace{1em}
+ c_{19} \tilde{p}\tilde{r}^2
+ c_{20} \tilde{r}^3. 
\end{split}
\Label{MS-h123}
\end{equation}
The unknown coefficients are found by substituting~(\ref{MS-h123}) in the invariance equation~(\ref{MS-invariance}) and equating coefficients of like monomials. 
One can show by induction on the power that quadratic and cubic terms dependent only on $\tilde{p}$ and $\tilde{r}$ must be zero. To show that there are no quadratic terms of this type, we examine the invariance equation (\ref{MS-invariance}) and observe that (i) $h_1 = 0$ implies that the function $n$ in the right member of the invariance equation cannot generate such a quadratic term, and (ii) the terms in the left member cannot generate such a quadratic term either. One can then show similarly
that there are no cubic terms of this type either. Hence, $b_8, b_9, b_{10} = 0$, and $c_{17}, c_{18}, c_{19}, c_{20} = 0$. The nonzero terms are
\begin{equation*}
\begin{split}
b_1 &= q (q+1)^3 s/ \ell^3, \\
b_2 &= 2 (q+1)^3 \left(q^3-q^2-q-1\right) s / \ell^4, \\
b_3 &= (q+1)^3 \left(q^6-2 q^5+2 q^3+4 q^2+2 q+1\right) s / q \ell^5 \\
b_4 &= q (q+1)^3 / \ell^3 , \\
b_5 &= (q+1)^3 \left(q^3-q^2-q-1 \right) / \ell^4 , \\
b_6 &= - b_4, \\
b_7 &= - q (q+1)^3 \left(2 q^2+q+1\right) / \ell^4 ,
\end{split}
\end{equation*}
and
\begin{equation*}
\begin{split}
c_1&= q (q+1)^3 / \ell^3 - 6q^2 (1+q)^6 s^2 / \ell^6 , \\
c_2&= (q-1) (q+1)^3 (4 q^2+3q+2) / \ell^4 - 2q(1+q)^6(9q^3-14q^2-14q-9)s^2 / \ell^7 , \\
c_3&= (q-1)^2 (q+1)^3 \left( 5 q^4+6q^3+5q^2+2q+1 \right) / q \ell^5 \\ & \qquad - 2q(1+q)^6 (9q^6-28q^5+2q^4+42q^3+58q^2+28q+9) s^2 / q \ell^8 , \\
c_4&= (q-1)^2 (q+1)^3 \left(2 q^6+q^5+q^4+3q^3+5q^2+3q+1\right) / q \ell^6 \\ & \qquad - 2 (q+1)^6 \left(3 q^9-14 q^8+16 q^7+21 q^6-32 q^5-92 q^4-81 q^3-44 q^2-14 q-3\right) s^2 / q \ell^9 , \\
c_5 &= - 9 q^2 (q+1)^6 s / \ell ^6, \\
c_6 &= - 2 q (q+1)^6 \left(9 q^3-14 q^2-14 q-9\right) s / \ell^7, \\
c_7 &= -(q+1)^6 \left(9 q^6-28 q^5+2 q^4+42 q^3+58 q^2+28 q+9\right) s / \ell^8 , \\
c_8 &= - q (q+1)^4 \left(q^4-5 q^3-13 q^2-5 q+1\right) s / \ell ^6 , \\
c_9 &= - 2 (q+1)^4 \left(q^7-9 q^6-17 q^5+4 q^4+16 q^3+9 q^2+q-1\right) s / \ell^7 , \\
c_{10} &= -(q+1)^4 \left(q^{10}-13 q^9-7 q^8+59 q^7+91 q^6+66 q^5+35 q^4+19 q^3+9 q^2+3 q+1\right) s / q \ell ^8, \\
c_{11}&= - 3 q^2 (q+1)^6 / \ell^6 , \\
c_{12}&= - q (q+1)^6 \left( 3q^3-4q^2-4q-3 \right) / \ell^7 , \\
c_{13}&= - q (q+1)^4 \left( q^4-2q^3-7q^2-2q+1 \right) / \ell^6 , \\
c_{14}&= - (q+1)^4 \left( q^7-6q^6-17q^5-11q^4-5q^3-3q^2-2q-1 \right) / \ell^7 , \\
c_{15}&= q (q+1)^4 \left( q^4+q^3-q^2+q+1\right) / \ell^6 , \\
c_{16}&= q (q+1)^4 \left( 2 q^6+3 q^5+3 q^3+3q^2+2q+1\right) / \ell^7.
\end{split}
\end{equation*}
These coefficients are used in the series representation~\eqref{MS-h} of $h(u,v,{\tilde p},{\tilde r})$ in Section~\ref{a4-ss-CenterManifolds}.

Finally, we observe that the MS model \eqref{MS-xyz'} is symmetric under the reflection
\[ (x,y,z,s) \mapsto (-x,-y,-z,-s), \]
and all terms in the center manifold expansion respect this symmetry. However, we recall that only the regime $s>0$ is of physical relevance in order to model the asymmetry between rapid deglaciation and slow glaciation.

\end{appendix}

\section*{Acknowledgments} 
The authors thank Edgar Knobloch for a useful conversation.

\small
\bibliography{DS_Paper}
\end{document}